\documentclass[11pt,a4paper]{article}

\usepackage[T1]{fontenc}
\usepackage[english]{babel}
\usepackage{amsmath,amsthm,amssymb,mathtools}
\usepackage{todonotes}

\usepackage[left=2.5cm,right=2.5cm,top=2.5cm,bottom=3cm]{geometry}
\usepackage[hidelinks]{hyperref}
\hypersetup{pdftitle={The small ball problem for the Brownian sheet: sharp bounds via entropy numbers in L-infinity}}
\numberwithin{equation}{section}
\newtheorem{theorem}{Theorem}[section]
\newtheorem{lemma}[theorem]{Lemma}
\newtheorem{corollary}[theorem]{Corollary}
\newtheorem{prop}[theorem]{Proposition}
\theoremstyle{remark}
\newtheorem{remark}[theorem]{Remark}
\newcommand{\E}{\mathbb E}
\newcommand{\Pp}{\mathbb P}
\newcommand{\R}{\mathbb R}
\newcommand{\N}{\mathbb N}
\newcommand{\one}{\mathbf 1}
\DeclareMathOperator{\Var}{Var}
\newcommand{\abs}[1]{\lvert#1\rvert}
\newcommand{\norm}[1]{\lVert#1\rVert}
\newcommand{\bW}{\mathbf{W}}
\newcommand{\bB}{\mathbf{B}}
\newcommand{\bH}{\mathbf{H}}
\newcommand{\D}{\mathcal D}
\DeclareMathOperator{\Hull}{Hull}

\title{The small ball problem for the Brownian sheet: sharp bounds via entropy numbers in $L_\infty$}
\author{Dmitriy Bilyk \and Sebastian Neumayer \and Tino Ullrich}
\date{\today}

\begin{document}
\maketitle

\begin{abstract}
We settle the uniform small deviation conjecture for the Brownian sheet in every dimension $d\ge2$, namely, that the logarithmic exponent is $2d-1$.
Our proof proceeds through sharp entropy and Kolmogorov estimates for the associated anchored mixed Sobolev embedding in the uniform norm.
Their connection to small deviations is known from the Kuelbs--Li theorem.
To establish the entropy estimates, the key is a probabilistic small ball inequality for smooth tensor product atoms whose exceptional set has Gaussian measure exponentially small in the number of atoms.
Its proof combines a count of dyadic rectangles with tools from probability, in particular Ehrhard's Gaussian Brunn--Minkowski inequality.
This also yields sharp entropy and Kolmogorov number estimates for embeddings of periodic mixed Sobolev, Besov and Nikol'skii classes with smoothness $r>1/2$ in $L_\infty$.
Finally, we extend our results to fractional and integrated Brownian sheets.  

\noindent\textbf{Lean certificate}: \url{https://github.com/Zeppo1994/brownian-sheet-certificate}
\end{abstract}

\section{Introduction and main results}\label{sec:intro}

The small deviation problem for the Brownian sheet asks for the behavior of the quantity
\begin{equation}\label{eq:smalldevfunction}
 \phi(\varepsilon)\coloneqq-\log\Pp\{\norm B_\infty\le\varepsilon\}
 \qquad\text{as }\varepsilon\downarrow0,
\end{equation}
where $B$ is the centered continuous Gaussian field on $[0,1]^d$ with covariance
\begin{equation}\label{eq:sheetcovariance}
 \E[B(\mathbf s)B(\mathbf t)]=\prod_{i=1}^d\min(s_i,t_i).
\end{equation}
For $d=1$ (Brownian motion), it is classical that $\phi(\varepsilon)\asymp\varepsilon^{-2}$.
Talagrand~\cite{Talagrand94} proved that $\phi(\varepsilon)\asymp\varepsilon^{-2}(\log(1/\varepsilon))^3$ for $d=2$.
For general $d\ge2$, the conjectured power of $\log(1/\varepsilon)$ is $2d-1$, and the matching upper estimate with this exponent is due to Dunker et al.~\cite{DKLL99}.
For $d\ge3$, the best lower estimate so far, by Bilyk, Lacey and Vagharshakyan~\cite[Thms.~2.10 and~2.13]{BLV08}, has exponent $2d-2+\eta(d)$ for some $0<\eta(d) < 1$ decaying with $d$.
For approaches to small deviations through chaining and path regularity, see Aurzada and Lifshits~\cite{AL08} and Aurzada~\cite{Aurzada11}.

Sharp lower estimates would follow from the small ball inequality for Haar functions, whose deterministic equivalent of \eqref{eq:SBI2} (see Remark~\ref{rem:SBI}) is, however, a central open problem, see \cite[Open Problem~1.2]{DTU18} and \cite[Conj.~4.4]{Bilyk11}.
In contrast, our approach is based on a probabilistic small ball inequality for hyperbolic sums of smooth atoms, valid for all coefficient vectors outside a set of exponentially small Gaussian measure (Theorem~\ref{lem:smoothconfinement}).
For $d=2$, an estimate of this type for random hyperbolic trigonometric polynomials follows from the deterministic small ball inequality, see~\cite[proof outline of Thm.~11.3.7]{TB04}.
In the cases where the deterministic inequality is open, we are unaware of such a result, neither for Haar functions nor for trigonometric polynomials. With the established probabilistic inequality, we can settle various open problems from approximation theory.
Similar random hyperbolic wavelet expansions have been used by Kempka, Schneider and Vyb\'iral~\cite{KSV24} to study the sample-path regularity of Brownian motion and the Brownian sheet in Besov-type spaces, including spaces of dominating mixed smoothness.
Horst, Jahn and Voigtlaender~\cite{HJV24} study the Besov regularity of random wavelet series on $\R^d$.

\begin{theorem}\label{thm:main}
As $\varepsilon\downarrow0$, the Brownian sheet in dimension $d\geq 2$ satisfies
\begin{equation}\label{eq:main}
 -\log\Pp\{\norm B_\infty\le\varepsilon\}
 \asymp_d \varepsilon^{-2}\bigl(\log(1/\varepsilon)\bigr)^{2d-1}.
\end{equation}
\end{theorem}

This result follows from the entropy lower bound in Theorem~\ref{thm:CM}, using the connection introduced by Kuelbs and Li~\cite{KL93} in the form established by Li and Linde~\cite[eq.~(1.3)]{LiLinde99}.
Section~\ref{sec:applications} derives further consequences for Gaussian fields and quantization.
In particular, Theorem~\ref{thm:fractionalsheet} gives sharp small-deviation orders for integrated Brownian sheets and fractional Brownian sheets whose smallest Hurst parameter occurs in at least two coordinates.
The fractional result matches the upper bound in~\cite[Cor.~7.5]{KuhnLinde02}.
The two-parameter fractional case was established by Belinsky and Linde~\cite{BelinskyLinde02}.
We also obtain sharp entropy and small-deviation orders for tensor-product Mat\'ern fields and sharp quantization and empirical approximation rates for the sheet laws.

\paragraph{Entropy and Kolmogorov numbers.}

Our proof strategy is connected to another well-known question in multivariate approximation theory.
Functions with dominating mixed smoothness are a classical model for high-dimensional approximation.
Their hyperbolic-cross structure keeps the main approximation rate independent of the dimension, and the dimension enters only through powers of the logarithm. 
Determining the powers is often the hard part, and it is hardest in the uniform norm.
There, sharp bounds for entropy and Kolmogorov numbers have been open for $d\ge3$, even for the simplest class with one mixed derivative in $L_2$.

The small-deviation-entropy theorem by Kuelbs and Li \cite{KL93} links the Brownian sheet to the $p=2$ case of the anchored integration operators $T_{d,p}\colon L_p([0,1]^d)\to C([0,1]^d)$, defined for $1<p<\infty$ by
\begin{equation}\label{eq:integration}
 (T_{d,p}g)(\mathbf t)\coloneqq\int_{[0,\mathbf t]}g(\mathbf u)\,d\mathbf u.
\end{equation}
We write $T_d\coloneqq T_{d,2}$.
Its compact image ball $K_d\coloneqq T_d(B_{L_2([0,1]^d)})$ is the Cameron--Martin unit ball of $B$ and can be written as the anchored mixed Sobolev unit ball
\begin{equation}\label{eq:Kd}
 K_d=\Bigl\{f\in C([0,1]^d):\ f(\mathbf t)=0\text{ if }t_i=0\text{ for some }i,\
 \norm{\partial_1\cdots\partial_d f}_{L_2}\le1\Bigr\},
\end{equation}
where $\partial_1\cdots\partial_d f$ denotes the weak mixed derivative.
The next theorem determines the entropy numbers $e_k$ of $T_{d,p}$ and its Kolmogorov numbers $d_k$, both as defined in~\eqref{eq:operatornumbers}.

\begin{theorem}\label{thm:CM}
Let $d\ge2$ and $1<p<\infty$.
The integration operator $T_{d,p}\colon L_p([0,1]^d)\to C([0,1]^d)$ defined by~\eqref{eq:integration} satisfies
\begin{equation}\label{eq:entropy}
 e_k(T_{d,p})\asymp_{d,p}k^{-1}(\log k)^{d-1/2}
 \qquad\text{as }k\to\infty.
\end{equation}
For $2\le p<\infty$, the same order holds for the Kolmogorov numbers $d_k(T_{d,p})$.
\end{theorem}
The upper bounds follow from Dunker et al.~\cite{DKLL99} for $p\ge2$ by inclusion and from Belinsky's periodic entropy estimate~\cite{Belinsky98} for $1<p<2$ by extension.
We establish the lower bounds in Section~\ref{sec:lower}.
The entropy estimate establishes Conjecture~2.9 in~\cite{BLV08} for $1<p<\infty$.

The same method applies to the periodic mixed Sobolev, Besov and H\"older-Nikol'skii spaces $\bW^r_p$, $\bB^r_{p,\theta}$ and $\bH^r_p=\bB^r_{p,\infty}$ on $\mathbb T^d$, as defined in Subsection~\ref{sec:spaces}, and their embedding into $L_\infty(\mathbb T^d)$.

\newpage
\begin{theorem}\label{thm:periodicintro}
Let $d\ge2$.
As $k\to\infty$, the following hold.
\begin{enumerate}
\item[(i)] For $2\leq p<\infty$ and $r>1/2$, it holds
\begin{equation}
    e_k(\operatorname I\colon\bW^r_p\to L_\infty)\asymp  d_k(\operatorname I\colon\bW^r_p\to L_\infty) \asymp k^{-r}(\log k)^{r(d-1)+1/2}.
\end{equation}
For $1 < p<\infty$ and $r>1/p$, the lower bound for $e_k$ remains valid.
\item[(ii)] For $2\le p\le\infty$, $2\le\theta\le\infty$ and $r>1/2$, it holds
\begin{equation}
    e_k(\operatorname I\colon\bB^r_{p,\theta}\to L_\infty)\asymp d_k(\operatorname I\colon\bB^r_{p,\theta}\to L_\infty)\asymp k^{-r}(\log k)^{(r+1/2-1/\theta)(d-1)+1/2}.
\end{equation}
The mixed H\"older--Nikol'skii spaces $\bH^r_p=\bB^r_{p,\infty}$ correspond to $\theta=\infty$.
For $1<p\le\infty$, $1\le\theta\le\infty$ and $r>1/p$, the lower bound for $e_k$ remains valid.
\end{enumerate}
The constants depend on $d,p,\theta,r$ and the choice of the equivalent norm.
\end{theorem}

The difficult part of Theorem~\ref{thm:periodicintro} is proving the lower bounds.
The mixed Sobolev and Nikol'skii upper bounds are known: for entropy numbers see Belinsky~\cite{Belinsky98} and~\cite[Thms.~11.2.4 and~11.2.10]{TB04}, and for Kolmogorov widths see~\cite[Thms.~11.3.4 and~11.3.5]{TB04}.
The two-dimensional lower bounds were also known, see~\cite[discussion after eq.~(1.3)]{Temlyakov17} and~\cite[Rem.~6.4]{TU22}.
For $d\ge3$, the two-sided estimates in Theorem~\ref{thm:periodicintro} settle the corresponding cases of~\cite[Open Problem~1.3]{DTU18}, \cite[Conj.~4.3.15]{DTU18} and the finite-$p$ cases of~\cite[Conj.~4.3.9]{DTU18}.
Moreover, part~(i) establishes the entropy order conjectured in~\cite[eq.~(1.8)]{Temlyakov17} and~\cite[Remark after Thm.~11.2.4]{TB04}, see Remark~\ref{rem:ZeroMean} for $1<p<2$.
In the small-smoothness range $2<p<\infty$, $1/p<r\le1/2$, the sharp order for mixed Sobolev classes remains open; for upper bounds, see~\cite[Thms.~6.1 and~6.3(i)]{TU22}.

For the mixed Besov spaces, real interpolation together with~\cite[Thm.~3.2]{TU22} gives the upper bounds in Theorem~\ref{thm:periodicintro}.
For $2\le p\le\infty$, $2\le\theta<\infty$ and $r>1/2$, Romanyuk~\cite[Thms.~3 and~3$'$]{Romanyuk19} proved the two-dimensional entropy and Kolmogorov-width orders.
Romanyuk and Yanchenko~\cite{RY23} obtained the corresponding entropy-number order for $d\ge2$ in the quasi-continuous norm.
Our uniform-norm estimate proves Romanyuk's Besov entropy conjecture~\cite[Conj.~1]{Romanyuk19} in this parameter range for $d\ge3$.
The range $\theta<1$, which is not covered here, is studied by Mayer and T.~Ullrich~\cite[Cor.~23, Rem.~24]{MU21}.
Further results on entropy numbers between nonperiodic mixed-smoothness function spaces can be found in Vyb\'iral~\cite{Vyb2006}.

The established width estimates also give bounds for sampling recovery.
Indeed, the Kolmogorov numbers of an operator coincide with the Kolmogorov widths of its image unit ball, and Temlyakov~\cite{Tem21} bounds the linear $L_2$ recovery error from $O(n)$ samples by the $n$-th uniform width, using weighted least squares.
Temlyakov and T.~Ullrich~\cite{TU21,TU22} apply this principle to classes with small mixed smoothness.
Further connections between approximation subspaces and recovery from function values are provided by sampling discretization~\cite{LT22} and the sampling projections of Krieg et al.~\cite{KPUU26}.
For a recent survey of optimal sampling and complexity, see Krieg and M.~Ullrich~\cite{KU26}.

For numerical integration, Novak~\cite{No86} bounds quadrature errors by Kolmogorov widths in $L_\infty$.
Going beyond this bound, Neumayer, Sch\"afer and T.~Ullrich~\cite{NeSchUl26} construct positive quadrature rules from random nodes and, after compression, obtain bounds governed by singular numbers of the Hilbert space embedding into $L_2$.
For periodic mixed Sobolev-Hilbert balls with $r>1/2$, their bound $n^{-r}(\log n)^{r(d-1)}$ is a factor $\sqrt{\log n}$ smaller than the sharp uniform-width order established in Theorem~\ref{thm:periodicintro}.
For these balls and every $d\ge2$, Novak's width-based bound therefore exceeds the optimal quadrature error by at least this factor.
For deterministic $L_2$ sampling recovery, Temlyakov's width-based bound exceeds the optimal error by exactly this factor~\cite{NaschUl22,DKU23,Tem21}.

\paragraph{Organization and notation.}

We introduce the dyadic atoms and the Gaussian field in Section~\ref{sec:sbi}, where we prove the probabilistic small ball inequality and discuss its deterministic counterpart.
A covering argument and Carl's inequality then yield Theorems~\ref{thm:CM} and~\ref{thm:periodicintro} in Section~\ref{sec:entropy}.
Identifying the Cameron--Martin space of the Brownian sheet allows us to deduce Theorem~\ref{thm:main} in Section~\ref{sec:sheet}.
The consequences for fractional integration and Brownian sheets, tensor-product Mat\'ern fields, and quantization are developed in Section~\ref{sec:applications}.

Throughout, $d\ge2$ is fixed, and $C,c>0$ may change from line to line and depend only on $\psi$ and $d$ unless stated otherwise.
All logarithms are natural.
We write $\abs{\mathbf s}\coloneqq \sum_i s_i$ for a multi-index $\mathbf s\in\N_0^d$, and $\lambda^d$ for Lebesgue measure on $[0,1]^d$.

\section{A probabilistic small ball inequality}\label{sec:sbi}\label{sec:atoms}\label{sec:psbi}
Here, we prove the probabilistic small ball inequality, representing the core of all the results in this paper.
The result covers any smooth atoms and does not assume moment conditions. 
Fix a nonzero real $\psi\in C_c^\infty(0,1)$ with $\norm\psi_\infty\le1$, and extend it by zero to $\R$.
For $s\in\N_0$ and $m=0,\ldots,2^s-1$, put
$I_{s,m}\coloneqq [m2^{-s},(m+1)2^{-s})$.
A dyadic rectangle is given as the tensor product  
\begin{equation}
R=R_{\mathbf s, \mathbf m} \coloneqq \prod_{i=1}^d I_{s_i,m_i} = \prod_{i=1}^d \bigl[m_i2^{-s_i},(m_i+1)2^{-s_i}\bigr)
\end{equation}
and has ``shape'' $\mathbf s = (s_1,...,s_d)$ and volume $\lambda^d(R)=2^{-\abs{\mathbf s}}$.
For $n\ge0$, let $\D_n$ be the set of  rectangles of volume $2^{-n}$.
The number of possible shapes is 
\begin{equation}\label{eq:shapes}
    S_n\coloneqq \binom{n+d-1}{d-1} \asymp_d n^{d-1}\,.
\end{equation}
Since there are $S_n$ shapes and $2^n$ rectangles of each shape, the number of rectangles in $\D_n$ satisfies $\vert \D_n \vert =2^nS_n \asymp 2^n n^{d-1}$.
For $R=R_{\mathbf s, \mathbf m}\in\D_n$, we define the atom
\begin{equation}\label{eq:atoms}
 a_R(\mathbf x)\coloneqq \prod_{i=1}^d\psi(2^{s_i}x_i-m_i).
\end{equation}
Each atom is smooth and satisfies $\operatorname{supp}a_R\subset R$ and
$\norm{a_R}_\infty\le1$.
Atoms of the same shape have disjoint supports.
Let $\mathbb T^d\coloneqq (\R/\mathbb Z)^d$ be equipped with the Lebesgue measure $\lambda^d$.
Since the atoms~\eqref{eq:atoms} vanish near the boundary of the cube $[0,1]^d$, we
may regard every finite sum of atoms as a smooth periodic function on
$\mathbb T^d$ with the same uniform norm.

Now, let $G_R$ with $R\in\D_n$ be independent standard Gaussian random variables on a probability space
$(\Omega,\mathcal F,\Pp)$, and put
\begin{equation}\label{eq:smoothlayer}
 V_n\coloneqq \sum_{R\in\D_n}G_Ra_R,
\end{equation}
with the atoms $a_R$ from~\eqref{eq:atoms}.
The following theorem is the probabilistic small ball inequality behind all lower bounds of this paper.
After discussing its connection to the deterministic case, we develop the required tools and then prove the theorem.

\begin{theorem}[Probabilistic small ball inequality]\label{lem:smoothconfinement}
There exist $a,b>0$, depending only on $d$ and $\psi$, such that
\begin{equation}\label{eq:smoothconfinement}
 \Pp\bigl\{\norm{V_n}_\infty\le a\sqrt{nS_n}\bigr\}\le\exp(-b\vert \D_n \vert)
\end{equation}
for all sufficiently large $n$.
\end{theorem}

\begin{remark}[Relation to the small ball inequality]\label{rem:SBI}
Consider the inequality
\begin{equation}\label{eq:SBI2}
 n^{(d-2)/2}\biggl\|\sum_{R\in\D_n}\lambda_Ra_R\biggr\|_\infty
 \ge c\,2^{-n}\sum_{R\in\D_n}\abs{\lambda_R}
 \qquad\text{for all }(\lambda_R)\in\R^{\vert \D_n \vert}.
\end{equation}
For tensor Haar functions, this is the fixed-volume form of the small ball conjecture~\cite[Open Problem~1.2]{DTU18}.
A smooth analogue with vanishing mean is formulated in~\cite[Conj.~2.5]{BLV08}. 
Both conjectures remain open for $d\ge3$. Besides small deviation probabilities for the Brownian sheet (hence the name) and the metric entropy of mixed smoothness classes, which we study here, such inequalities play a central role in discrepancy
theory~\cite{Bilyk11,Bilyk14,BilykLacey13}. 
Let $E_n$ be the set of coefficient vectors violating~\eqref{eq:SBI2}.
For $\lambda=\mathbf G$, the right-hand side is at most $2cS_n$ outside an event of probability $e^{-c'\vert \D_n \vert}$, by Gaussian concentration of $\mathbf G\mapsto\norm{\mathbf G}_1$.
Since $n^{-(d-2)/2}S_n\asymp\sqrt{nS_n}$, for small $c$ every $\mathbf G\in E_n$ outside this event satisfies $\norm{V_n}_\infty\le a\sqrt{nS_n}$.
A union bound and Theorem~\ref{lem:smoothconfinement} give
\begin{equation}
 \Pp\{\mathbf G\in E_n\}\le e^{-c'\vert \D_n \vert}+e^{-b\vert \D_n \vert}\le e^{-c''\vert \D_n \vert}.
\end{equation}
Thus~\eqref{eq:SBI2} holds for all coefficient vectors outside a set of
Gaussian measure $e^{-c''\vert \D_n \vert}$.
\end{remark}

\paragraph{Proof strategy.}

The proof of Theorem~\ref{lem:smoothconfinement} relies on the following three ingredients:

\begin{enumerate}
\item[(i)] \emph{Good regions.}
At each level $n$, we split the cube $[0,1]^d$ according to whether the variance of the Gaussian field $V_n$ is at least half its spatial mean or below it.
The first region, denoted by $A_n$, occupies a fixed positive fraction of the cube.
On $A_n$, we divide the field $V_n$ by its pointwise standard deviation.

\item[(ii)] \emph{Average the peaks.}
We form a symmetric exponential average $U_\tau$ of the normalized Gaussian field on $A_n$, with $\E[U_\tau]=1$.
Then, we write its variance as a double average of $\cosh(\tau\rho)-1$, where $\rho$ denotes correlation.
Grouping pairs of points by their dyadic hulls controls both their correlations (Lemma~\ref{lem:correlations}) and the measure of the corresponding pairs.
Since $\lambda^d(A_n)$ is bounded below, these estimates yield $\Var(U_\tau)\lesssim\tau^2/\lvert\D_n\rvert$ for $0<\tau\le cn$ (Lemma~\ref{lem:expvariance}).

\item[(iii)] \emph{Gaussian convexity.}
We now choose $\tau \asymp n$ according to (ii).
Then $U_\tau$ concentrates near its mean one as $n\to\infty$.
On the small-ball event $\{\norm{V_n}_\infty\le a\sqrt{nS_n}\}$, however, $U_\tau\le e^{-\tau/4}\to0$.
Since $\log U_\tau$ is convex, Proposition~\ref{lem:convextail} bounds the probability of this downward deviation by $\exp(-c\tau^2/\Var(U_\tau))\le e^{-c'\lvert\D_n\rvert}$, using the estimate from (ii).
\end{enumerate}

\subsection{Correlations and dyadic hulls}

First, we define the variance $v_n(\mathbf x)\coloneqq \E[V_n(\mathbf x)^2]$ and its spatial mean \smash{$\bar v_n\coloneqq \int_{[0,1]^d}v_n\,d\lambda^d$}.
At each point $\mathbf x \in [0,1)^d$, at most one atom from each of the $S_n$ shapes is nonzero.
Thus, independence of the Gaussian coefficients and $\norm\psi_\infty\le1$ give
\begin{equation}\label{eq:variance}
 v_n(\mathbf x)=\sum_{R\in\D_n}a_R(\mathbf x)^2\le S_n,
 \qquad \bar v_n=\norm{\psi}_2^{2d}S_n.
\end{equation}
Further, we define our candidate set for peaks of $V_n$ and an associated restricted measure as
\begin{equation}\label{eq:goodset}
 A_n\coloneqq \bigl\{\mathbf x\in[0,1]^d:v_n(\mathbf x)\ge\bar v_n/2\bigr\},\qquad
 \nu_n(E)\coloneqq \frac{\lambda^d(E\cap A_n)}{\lambda^d(A_n)}
 \quad \forall \, E\subseteq[0,1]^d\text{ measurable}.
\end{equation}
Set $c_{d,\psi}\coloneqq\frac12\norm{\psi}_2^{2d}\in(0,1/2]$. Then, the inequality $\bar v_n\le\bar v_n/2+\lambda^d(A_n)\norm{v_n}_\infty$ and \eqref{eq:variance} imply
\begin{equation}\label{eq:goodvariance}
 \lambda^d(A_n)\ge c_{d,\psi},\qquad v_n(\mathbf x)\ge c_{d,\psi}S_n
 \qquad \forall\,\mathbf x\in A_n.
\end{equation}
For $\mathbf x\in A_n$, the normalized variable $Y_{\mathbf x}\coloneqq V_n(\mathbf x)/\sqrt{v_n(\mathbf x)}$ is standard Gaussian.
Moreover, each $Y_{\mathbf x}$ is a linear functional of $\mathbf G=(G_R)_{R\in\D_n}$.
By independence, two such variables have correlation
\begin{equation}\label{eq:rho}
 \rho(\mathbf x,\mathbf y)\coloneqq\E[Y_{\mathbf x}Y_{\mathbf y}]
 =\frac{1}{\sqrt{v_n(\mathbf x)v_n(\mathbf y)}}\sum_{R\in\D_n}a_R(\mathbf x)a_R(\mathbf y).
\end{equation}
A summand in~\eqref{eq:rho} can be nonzero only if $R$ contains both $\mathbf x$ and $\mathbf y$.

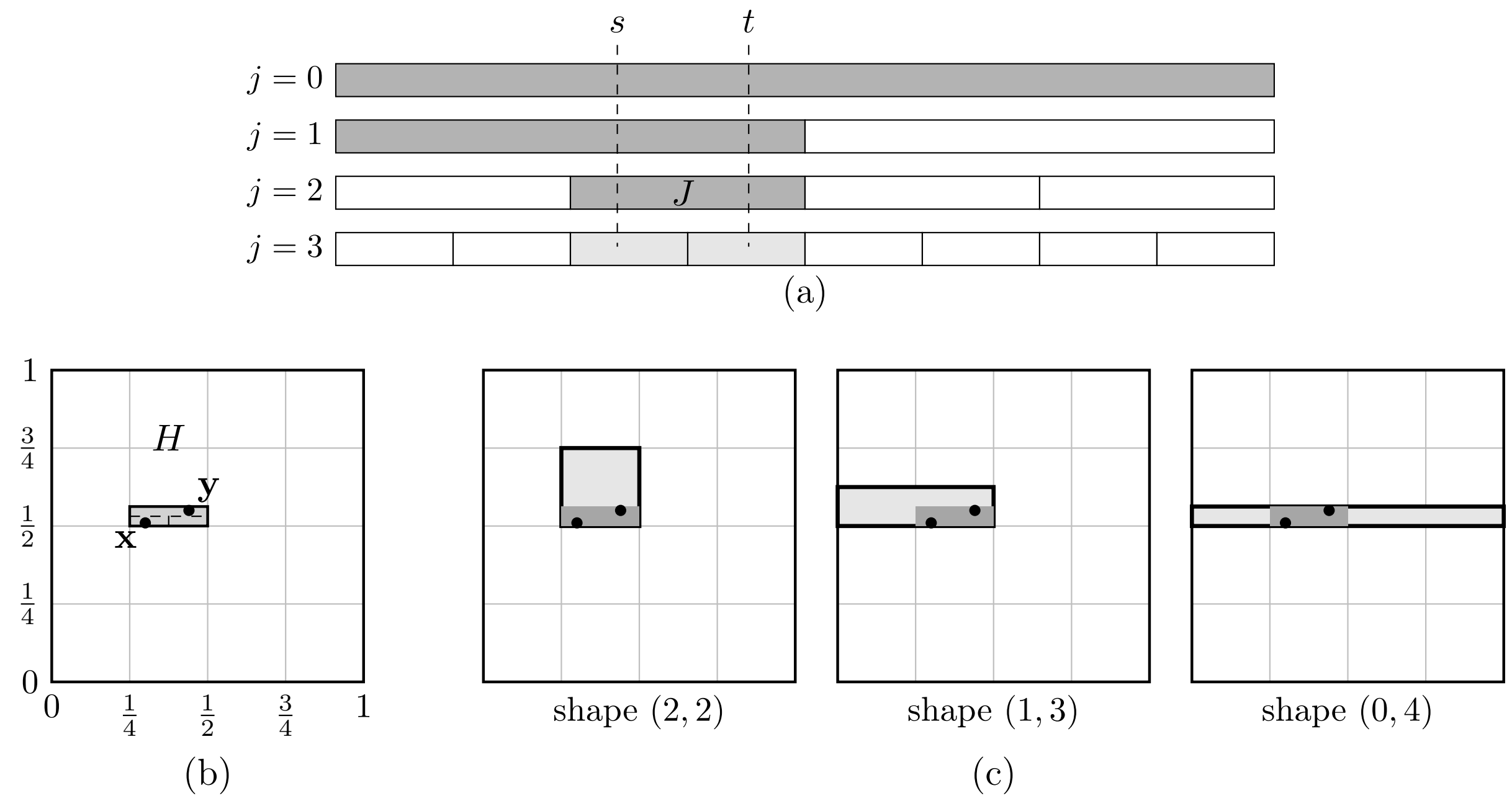
\begin{figure}[t]
\centering
\begin{tikzpicture}[xscale=10,yscale=1]
 \fill[gray!60] (0,0) rectangle (1,0.35);
 \fill[gray!60] (0,-0.6) rectangle (0.5,-0.25);
 \fill[gray!60] (0.25,-1.2) rectangle (0.5,-0.85);
 \fill[gray!20] (0.25,-1.8) rectangle (0.5,-1.45);
 \foreach \j in {0,...,3}{
   \pgfmathsetmacro{\w}{1/2^\j}
   \pgfmathtruncatemacro{\N}{2^\j-1}
   \foreach \m in {0,...,\N}{
     \draw (\m*\w,-\j*0.6) rectangle ({(\m+1)*\w},-\j*0.6+0.35);}
   \node[left] at (0,-\j*0.6+0.175) {\small $j=\j$};}
 \draw[dashed] (0.30,0.55) -- (0.30,-1.6);
 \draw[dashed] (0.44,0.55) -- (0.44,-1.6);
 \node[above] at (0.30,0.55) {$s$};
 \node[above] at (0.44,0.55) {$t$};
 \node at (0.37,-1.025) {$J$};
 \node at (0.5,-2.1) {(a)};
\end{tikzpicture}

\medskip

\resizebox{\linewidth}{!}{%
\begin{tikzpicture}[scale=3.2,baseline=0pt]
 \foreach \t in {0.25,0.5,0.75}{
   \draw[gray!50,thin] (\t,0) -- (\t,1);
   \draw[gray!50,thin] (0,\t) -- (1,\t);}
 \fill[gray!35] (0.25,0.5) rectangle (0.5,0.5625);
 \draw[thick] (0,0) rectangle (1,1);
 \draw[thick] (0.25,0.5) rectangle (0.5,0.5625);
 \draw[dashed] (0.375,0.5) -- (0.375,0.5625);
 \draw[dashed] (0.25,0.53125) -- (0.5,0.53125);
 \fill (0.30,0.51) circle (0.018);
 \fill (0.44,0.55) circle (0.018);
 \node[below left=-2pt] at (0.30,0.51) {$\mathbf x$};
 \node[above right=-2pt] at (0.44,0.55) {$\mathbf y$};
 \node[above] at (0.375,0.70) {$H$};
 \foreach \t/\l in {0/0,0.25/\tfrac14,0.5/\tfrac12,0.75/\tfrac34,1/1}{
   \node[below] at (\t,0) {\small$\l$};
   \node[left] at (0,\t) {\small$\l$};}
 \node at (0.5,-0.30) {(b)};
\end{tikzpicture}%
\hspace{1cm}%
\begin{tikzpicture}[scale=3.2,baseline=0pt]
 \foreach \t in {0.25,0.5,0.75}{
   \draw[gray!50,thin] (\t,0) -- (\t,1);
   \draw[gray!50,thin] (0,\t) -- (1,\t);}
 \fill[gray!20] (0.25,0.5) rectangle (0.5,0.75);
 \draw[very thick] (0.25,0.5) rectangle (0.5,0.75);
 \fill[gray!70] (0.25,0.5) rectangle (0.5,0.5625);
 \draw[thick] (0,0) rectangle (1,1);
 \fill (0.30,0.51) circle (0.018);
 \fill (0.44,0.55) circle (0.018);
 \node[below] at (0.5,0) {\small shape $(2,2)$};
\end{tikzpicture}%
\hspace{0.4cm}%
\begin{tikzpicture}[scale=3.2,baseline=0pt]
 \foreach \t in {0.25,0.5,0.75}{
   \draw[gray!50,thin] (\t,0) -- (\t,1);
   \draw[gray!50,thin] (0,\t) -- (1,\t);}
 \fill[gray!20] (0,0.5) rectangle (0.5,0.625);
 \draw[very thick] (0,0.5) rectangle (0.5,0.625);
 \fill[gray!70] (0.25,0.5) rectangle (0.5,0.5625);
 \draw[thick] (0,0) rectangle (1,1);
 \fill (0.30,0.51) circle (0.018);
 \fill (0.44,0.55) circle (0.018);
 \node[below] at (0.5,0) {\small shape $(1,3)$};
 \node at (0.5,-0.30) {(c)};
\end{tikzpicture}%
\hspace{0.4cm}%
\begin{tikzpicture}[scale=3.2,baseline=0pt]
 \foreach \t in {0.25,0.5,0.75}{
   \draw[gray!50,thin] (\t,0) -- (\t,1);
   \draw[gray!50,thin] (0,\t) -- (1,\t);}
 \fill[gray!20] (0,0.5) rectangle (1,0.5625);
 \draw[very thick] (0,0.5) rectangle (1,0.5625);
 \fill[gray!70] (0.25,0.5) rectangle (0.5,0.5625);
 \draw[thick] (0,0) rectangle (1,1);
 \fill (0.30,0.51) circle (0.018);
 \fill (0.44,0.55) circle (0.018);
 \node[below] at (0.5,0) {\small shape $(0,4)$};
\end{tikzpicture}%
}
\caption{Dyadic chains and hulls.
(a) The intervals containing both points form a chain ending in $J$.
(b) The dyadic hull $H = \mathrm{Hull}(\mathbf x, \mathbf y)$ of the two points $\mathbf x$ and $\mathbf y$.
(c) The three dyadic rectangles of volume $2^{-4}$ containing $H$, with their shapes indicated below.}
\label{fig:chain}
\label{fig:hull}
\end{figure}

Recall that any two dyadic intervals are either nested or disjoint.
For two points $s\ne t$ in $[0,1)$, the dyadic intervals containing $s$ and $t$ form a chain ending in a smallest interval $J$, whose halves separate $s$ and $t$, see Figure~\ref{fig:chain}(a).
For $\mathbf x,\mathbf y\in[0,1)^d$ with $x_i\ne y_i$ for every $i$, let $J_i$ be the smallest dyadic interval containing $x_i$ and $y_i$, and define
\begin{equation}\label{eq:hull}
 \Hull(\mathbf x,\mathbf y)\coloneqq J_1\times\cdots\times J_d.
\end{equation}
This is the smallest dyadic rectangle containing both points, see Figure~\ref{fig:hull}(b).
Since small hulls permit many shared atoms and hence large correlation, we estimate how rarely such pairs occur.

\begin{lemma}\label{lem:correlations}
For every $u \in \N_0$, it holds
\begin{equation}\label{eq:prefixlaw}
 \lambda^{2d}\{\Hull(\mathbf x,\mathbf y)\in\D_u\}
      =2^{-u-d}S_u.
\end{equation}
Moreover, for almost every pair $(\mathbf x,\mathbf y)\in A_n^2$, it holds
$\rho(\mathbf x,\mathbf y)=0$ if $\lambda^d(\Hull(\mathbf x,\mathbf y))>2^{-n}$, and
\begin{equation}\label{eq:correlation}
 |\rho(\mathbf x,\mathbf y)|\le\min\{1,c_{d,\psi}^{-1}S_h/S_n\}
 \qquad\text{if }\Hull(\mathbf x,\mathbf y)\in\D_{n+h},\ h\ge0.
\end{equation}
\end{lemma}

\begin{proof}
For a dyadic rectangle $H=\prod_iJ_i\subset[0,1)^d$, the dyadic hull of $(\mathbf x, \mathbf y) \in [0,1)^{2d}$ is $H$ exactly when $x_i$ and $y_i$ lie in opposite halves of $J_i$ for every $i$.
Each coordinate falls in a prescribed half with probability $|J_i|/2$, so the two possible orders in each coordinate give
\begin{equation}
 \lambda^{2d}\{\Hull(\mathbf x,\mathbf y)=H\}
 =\prod_{i=1}^d 2\left(\frac{|J_i|}{2}\right)^2
 =2^{-d}\lambda^d(H)^2.
\end{equation}
Taking the union over the $\vert \D_u \vert=2^uS_u$ rectangles $H\in\D_u$, each of volume $2^{-u}$, gives \eqref{eq:prefixlaw}.

For $(\mathbf x,\mathbf y)\in A_n^2$ with $x_i\ne y_i$ for every $i$, put $H=\Hull(\mathbf x,\mathbf y)$.
A rectangle $R\in\D_n$ contains both points if and only if $R\supseteq H$.
If $\lambda^d(H)>2^{-n}$, there is no such rectangle, so~\eqref{eq:rho} gives $\rho(\mathbf x,\mathbf y)=0$.
If $\lambda^d(H)=2^{-n-h}$ with $h\ge0$, each such rectangle is determined by $k_1,\ldots,k_d\in\N_0$ doublings of the sides of $H$, with $k_1+\cdots+k_d=h$, see Figure~\ref{fig:hull}(c).
There are at most $\binom{h+d-1}{d-1}=S_h$ such rectangles.
Since $\norm{a_R}_\infty\le1$, equations~\eqref{eq:rho} and~\eqref{eq:goodvariance} give
\begin{equation}
 \abs{\rho(\mathbf x,\mathbf y)}
 \le\frac{\#\{R\in\D_n:R\supseteq H\}}{c_{d,\psi}S_n}
 \le\frac{S_h}{c_{d,\psi}S_n}.
\end{equation}
Together with $\abs\rho\le1$, this proves~\eqref{eq:correlation}.
\end{proof}

\subsection{The soft maximum and its variance}
To bound the small-ball probability, we now construct a random variable that becomes small whenever $V_n$ is uniformly small.
The hull and correlation estimates \eqref{eq:prefixlaw} and \eqref{eq:correlation} will show that this variable is typically close to one, making such a deviation unlikely.
For $\tau>0$, we define
\begin{equation}\label{eq:expaverage}
 U_\tau\coloneqq e^{-\tau/2}
 \int_{A_n}\cosh(\sqrt\tau\,Y_{\mathbf x})\,d\nu_n(\mathbf x).
\end{equation}
Writing $\cosh u=(e^u+e^{-u})/2$, we see that
$\tau^{-1/2}\log U_\tau$ is, up to an additive constant, a soft maximum of the family
$\{Y_{\mathbf x},-Y_{\mathbf x}:\mathbf x\in A_n\}$.
In the covariance, this symmetrization cancels the term linear in $\rho$, yielding the variance bound in the following lemma.

\begin{lemma}\label{lem:expvariance}
For every $\tau>0$, it holds $\E[U_\tau]=1$, and $\log U_\tau$ is a finite convex function of $\mathbf G$.
Moreover, $\E[U_\tau^p]<\infty$ for every $p>0$.
Further, there exist $\kappa_{d,\psi}, C_{d,\psi}>0$ such that for every $0<\tau\le\kappa_{d,\psi}(n+1)$ it holds
\begin{equation}\label{eq:expvariance}
 \Var(U_\tau)\le C_{d,\psi}\frac{\tau^2}{\vert \D_n \vert}.
\end{equation}
\end{lemma}

\begin{proof}
The Gaussian exponential formula, namely $\E[e^{tZ}]=e^{t^2/2}$ for $Z\sim N(0,1)$ and $t \in \R$, gives $\E[\cosh(\sqrt\tau\,Y_{\mathbf x})]=e^{\tau/2}$ and hence $\E[U_\tau]=1$.
Recall that $Y_{\mathbf x}(\mathbf G)= v_n(\mathbf x)^{-1/2}\sum_{R\in\D_n}G_Ra_R(\mathbf x)$ is linear in the coefficients $\mathbf G$.
Writing $\cosh u=(e^u+e^{-u})/2$ in~\eqref{eq:expaverage} and applying H\"older's inequality gives
\begin{equation}
 U_\tau\bigl(\theta\mathbf G_1+(1-\theta)\mathbf G_2\bigr)
 \le U_\tau(\mathbf G_1)^\theta
      U_\tau(\mathbf G_2)^{1-\theta}
\end{equation}
for arbitrary vectors $\mathbf G_1,\mathbf G_2\in\R^{\lvert\D_n\rvert}$ and $0<\theta<1$.
Taking logarithms shows that $\log U_\tau$ is convex as a function of the coefficients $\mathbf G$.
Moreover, by Cauchy--Schwarz, it holds $\abs{Y_{\mathbf x}}\le\norm{\mathbf G}_2$.
Since $1\le\cosh u\le e^{\abs u}$ and $\nu_n(A_n)=1$, we thus obtain that
\begin{equation}
 -\tau/2\le\log U_\tau
 \le-\tau/2+\sqrt\tau\,\norm{\mathbf G}_2.
\end{equation}
This proves that $\log U_\tau$ is finite and that $\E[U_\tau^p] \le e^{-p\tau/2}\E [e^{p\sqrt\tau\,\|\mathbf G\|_2}] <\infty$ for every $p>0$.

Since $\E[U_\tau]=1$, applying Tonelli's theorem to the nonnegative integrand in~\eqref{eq:expaverage} gives
\begin{equation}\label{eq:VarianceFormula}
 \Var(U_\tau) =\E[U_\tau^2]-1 =e^{-\tau}\int_{A_n^2}
 \E\bigl[\cosh(\sqrt\tau\,Y_{\mathbf x})
          \cosh(\sqrt\tau\,Y_{\mathbf y})\bigr] \,d\nu_n(\mathbf x)\,d\nu_n(\mathbf y) - 1.
\end{equation}
For $\mathbf x,\mathbf y\in A_n$, the jointly Gaussian variables $Y_{\mathbf x},Y_{\mathbf y}\sim\mathcal N(0,1)$ satisfy the identity $\Var(Y_{\mathbf x}\pm Y_{\mathbf y})=2(1\pm\rho(\mathbf x,\mathbf y))$.
Using $\cosh a\cosh b=(\cosh(a+b)+\cosh(a-b))/2$, the Gaussian exponential formula thus gives
\begin{align}
\E\bigl[\cosh(\sqrt\tau\,Y_{\mathbf x})\cosh(\sqrt\tau\,Y_{\mathbf y})\bigr] &=\frac{1}{4}\E\Bigl[
 e^{\sqrt\tau(Y_{\mathbf x}+Y_{\mathbf y})}+e^{-\sqrt\tau(Y_{\mathbf x}+Y_{\mathbf y})}
 +e^{\sqrt\tau(Y_{\mathbf x}-Y_{\mathbf y})}+e^{-\sqrt\tau(Y_{\mathbf x}-Y_{\mathbf y})} \Bigr]\notag\\
 &=\frac{1}{2}
 \bigl(e^{\tau(1+\rho(\mathbf x,\mathbf y))}+e^{\tau(1-\rho(\mathbf x,\mathbf y))}\bigr)
 =e^{\tau}\cosh(\tau\rho(\mathbf x,\mathbf y)).
 \label{eq:expcovariance}
\end{align}
Inserting this into \eqref{eq:VarianceFormula} gives
\begin{equation}\label{eq:IntegralVariance}
    \Var(U_\tau)=\int_{A_n^2}\bigl(\cosh(\tau\rho(\mathbf x,\mathbf y))-1\bigr)
 \,d\nu_n(\mathbf x)\,d\nu_n(\mathbf y).
\end{equation}
Now, for $h\ge0$, put $E_h\coloneqq\{(\mathbf x,\mathbf y)\in A_n^2: \Hull(\mathbf x,\mathbf y)\in\D_{n+h}\}$.
By Lemma~\ref{lem:correlations}, $\rho=0$ when the hull has volume greater than $2^{-n}$.
Thus, grouping the integral in \eqref{eq:IntegralVariance} over the sets $E_h$ and using that $\cosh u-1\le u^2e^{\abs u}/2$ gives
\begin{equation}\label{eq:GroupedIntegralEstimate}
 \Var(U_\tau)\le\frac{\tau^2}{2}\sum_{h\ge0}
 \int_{E_h}\rho(\mathbf x,\mathbf y)^2e^{\tau\abs{\rho(\mathbf x,\mathbf y)}}
 \,d\nu_n(\mathbf x)\,d\nu_n(\mathbf y).
\end{equation}

By \eqref{eq:goodset} and \eqref{eq:goodvariance}, $\nu_n$ has density at most $c_{d,\psi}^{-1}$ with respect to $\lambda^d$.
Thus $\nu_n\otimes\nu_n$ has density at most $c_{d,\psi}^{-2}$ with respect to $\lambda^{2d}$, and the hull estimate~\eqref{eq:prefixlaw} yields
\begin{equation}\label{eq:Ehmeasure}
 (\nu_n\otimes\nu_n)(E_h)\le c_{d,\psi}^{-2}2^{-n-h}S_{n+h}.
\end{equation}
The correlation bound~\eqref{eq:correlation} gives $\abs\rho\le c_{d,\psi}^{-1}S_h/S_n$ on $E_h$.
For $h\le n$, we have by \eqref{eq:shapes} that
\begin{equation}\label{eq:ProdFormula}
 \frac{S_h}{S_n}
 =\prod_{i=1}^{d-1}\frac{h+i}{n+i}
 \le\frac{h+1}{n+1},
\end{equation}
which implies $\abs\rho\le c_{d,\psi}^{-1}(h+1)/(n+1).$
For $h>n$, the same bound follows from $\abs\rho\le1$ and $c_{d,\psi}^{-1}\ge2$.
Thus, setting $\kappa_{d,\psi}\coloneqq\frac12c_{d,\psi}\log2$ leads to $e^{\tau\abs\rho}\le2^{(h+1)/2}$ for $0<\tau\le\kappa_{d,\psi}(n+1)$.

Substituting \eqref{eq:Ehmeasure} together with the bounds on $\rho$ and $e^{\tau\abs\rho}$ into \eqref{eq:GroupedIntegralEstimate}, and using $2^{1/2}/2\le1$, gives
\begin{equation}\label{eq:variancesum}
 \Var(U_\tau)
 \le c_{d,\psi}^{-4}\tau^2\sum_{h\ge0}2^{-n-h}S_{n+h}(S_h/S_n)^2\,2^{h/2}.
\end{equation}
Since $(n+h+i)i\le(n+i)(h+i)$ for $i=1,\ldots,d-1$, we can multiply these inequalities to obtain $S_{n+h}\le S_nS_h$.
Inserting this estimate and $\vert\D_n\vert=2^nS_n$ into~\eqref{eq:variancesum}, we obtain
\begin{equation}
 \Var(U_\tau)
 \le\frac{c_{d,\psi}^{-4}\tau^2}{\vert\D_n\vert}\sum_{h\ge0}2^{-h/2}S_h^3.
\end{equation}
Since $S_h\le(h+1)^{d-1}$, see \eqref{eq:shapes}, the series on the right converges.
We now set $C_{d,\psi}\coloneqq c_{d,\psi}^{-4}\sum_{h\ge0}2^{-h/2}S_h^3$ and obtain $\Var(U_\tau)\le C_{d,\psi}\frac{\tau^2}{\vert\D_n\vert}$ as desired.
\end{proof}

\subsection{A Gaussian lower tail bound for convex functions}\label{sec:convextail}

This subsection is independent of the small ball problem.
We show that a convex function $f$ of a Gaussian vector has a Gaussian lower tail when $\E[e^{f(\mathbf G)}]=1$.
The tail bound depends only on the second moment of $e^{f(\mathbf G)}$.
The underlying tool is Ehrhard's Gaussian Brunn--Minkowski inequality~\cite{Ehrhard83}.
Let $\Phi$ denote the standard Gaussian distribution function with $\Phi^{-1}(0)=-\infty$ and $\Phi^{-1}(1)=\infty$.

\begin{lemma}[Ehrhard]\label{lem:ehrhard}
Let $\gamma_D$ be the standard Gaussian measure on $\R^D$, and let $K,L\subset\R^D$ be closed convex sets with $0<\gamma_D(K),\gamma_D(L)<1$.
Then, for every $\lambda\in[0,1]$, it holds
\begin{equation}\label{eq:ehrhard}
 \Phi^{-1}\bigl(\gamma_D(\lambda K+(1-\lambda)L)\bigr)
 \ge\lambda\Phi^{-1}(\gamma_D(K))
 +(1-\lambda)\Phi^{-1}(\gamma_D(L)).
\end{equation}
\end{lemma}

The following proposition combines Ehrhard's inequality with Chebyshev's inequality to obtain a Gaussian lower-tail bound.

\begin{prop}[Lower tail of a convex logarithm]\label{lem:convextail}
Let $\mathbf G$ be a standard Gaussian vector in $\R^D$, and let $f\colon\R^D\to\R$ be convex.
Suppose that $\E[e^{f(\mathbf G)}]=1$ and $\Var(e^{f(\mathbf G)})\le\sigma^2$ for some $0<\sigma\le1/4$.
Then, with $q\coloneqq\Phi^{-1}(3/4)$, it holds
\begin{equation}\label{eq:convextail}
 \Pp\{f(\mathbf G)\le-t\}
 \le\Phi\Bigl(-\frac{q(t-1)}{3\sigma}\Bigr)
 \le\exp\Bigl(-\frac{q^2(t-1)^2}{18\sigma^2}\Bigr)
 \qquad\forall \, t\ge1.
\end{equation}
\end{prop}

\begin{figure}[t]
\centering
\begin{tikzpicture}[xscale=2.2,yscale=0.42,
  declare function={G(\s)=2*\s-0.6*(\s*\s-0.25);}]
 \draw[->] (-2.9,0) -- (1.3,0) node[right] {$s$};
 \draw[->] (0,-9.3) -- (0,3.4);

 \draw[gray,dashed] (0,1) -- (1.2,1)
   node[right,black] {\small$\widetilde F(s_2)\ge q$};
 \draw[gray,dashed] (0,-1) -- (-1.1,-1)
   node[left,black] {\small$\widetilde F(s_1)\le-q$};

 \draw[thick,dashed] (-2.5,-5) -- (1.0,2);
 \draw[very thick,domain=-2.55:1.05,samples=80]
   plot (\x,{G(\x)});

 \fill (-0.5,-1) circle (0.04 and 0.21);
 \fill (0.5,1) circle (0.04 and 0.21);
 \draw (-0.5,0.2) -- (-0.5,-0.2);
 \draw (0.5,0.2) -- (0.5,-0.2);
 \node[above] at (-0.5,0.2) {\small$s_1$};
 \node[below right=-1pt] at (0.5,-0.05) {\small$s_2$};
 \draw[densely dotted] (-0.5,-1) -- (-0.5,-2.4);
 \draw[densely dotted] (0.5,1) -- (0.5,-2.4);
 \draw[<->] (-0.5,-2.4) -- (0.5,-2.4)
   node[right] {\small$\le6\sigma$};

 \draw (-2.5,0.2) -- (-2.5,-0.2);
 \node[above] at (-2.5,0.2) {\small$-t$};
 \draw[densely dotted] (-2.5,0) -- (-2.5,{G(-2.5)});
 \fill (-2.5,-5) circle (0.04 and 0.21);
 \fill (-2.5,{G(-2.5)}) circle (0.04 and 0.21);
 \node[left] at (-2.5,-5) {\small secant};
 \node[left] at (-2.5,{G(-2.5)})
   {\small$\widetilde F(-t)$};
 \node[above left] at (-1.6,-3.4)
  {\small slope $\ge q/(3\sigma)$};
\end{tikzpicture}
\caption{Schematic illustration of the secant argument in
Proposition~\ref{lem:convextail}.
To the left of $s_1$, concavity places the graph below the dashed secant extension, yielding $\widetilde F(-t)\le-q(t-1)/(3\sigma)$.}
\label{fig:secant}
\end{figure}
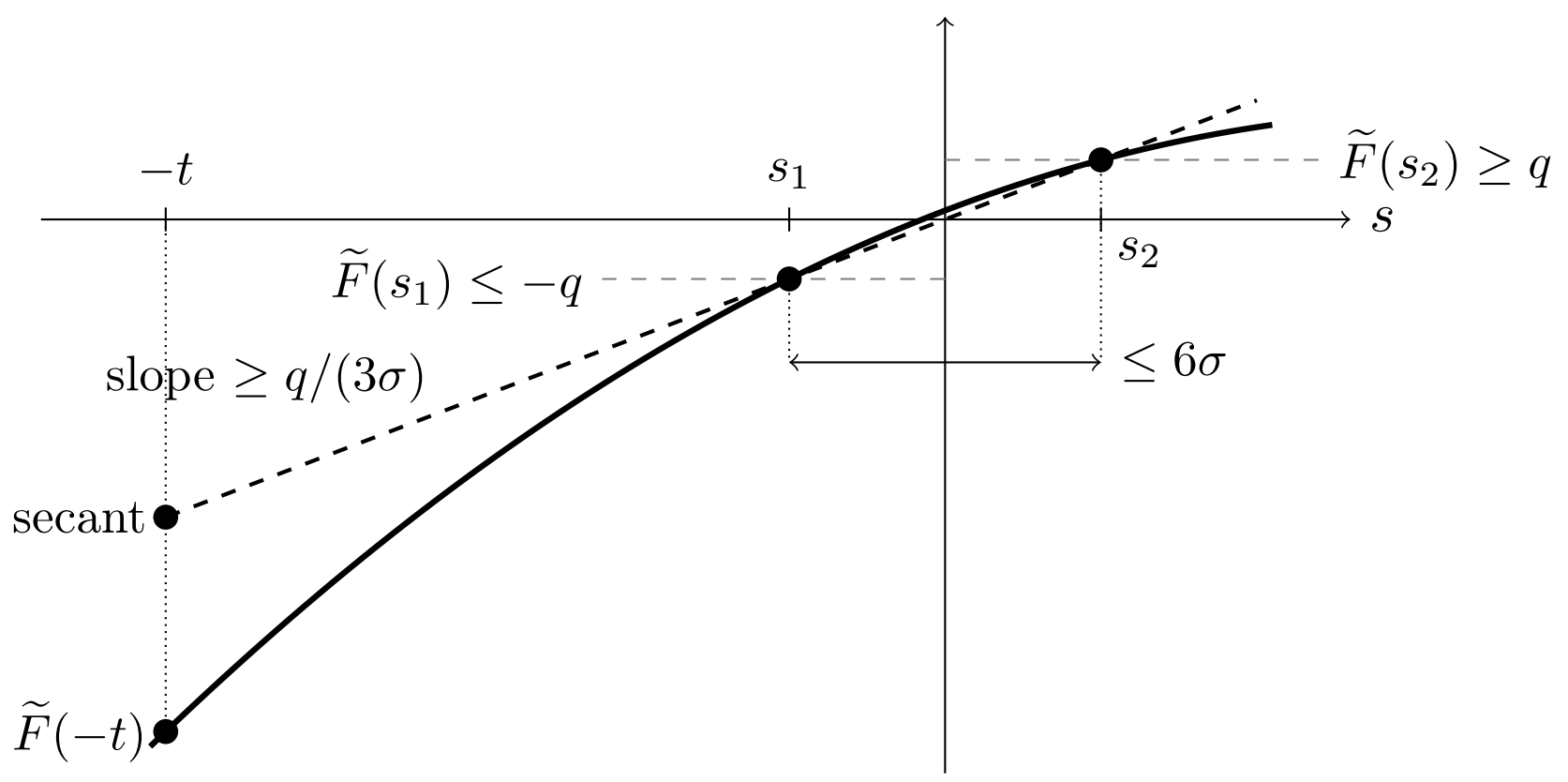

\begin{proof}
Fix $t\ge1$ and define $F(s)\coloneqq\Pp\{f(\mathbf G)\le s\}$.
If $F(-t)=0$, there is nothing to prove.
Otherwise, the normalization $\E[e^{f(\mathbf G)}]=1$ implies that $f$ is nonconstant.
Then, by finiteness and convexity, $f$ is continuous and unbounded above.
In particular, we get $F(s)<1$ for $s\in\R$.

Below, we will bound $F(-t)$ by a secant argument, see Figure~\ref{fig:secant} for an illustration.
To this end, we first introduce $\widetilde F\coloneqq\Phi^{-1}\circ F$ and establish its concavity.
Since $0<F(-t)\le F(s)<1$ for $s\ge-t$, this function is finite on $[-t,\infty)$.
By continuity and convexity of $f$, the sublevel sets $K_s\coloneqq\{f\le s\}$ are closed and convex.
For $s,s'\ge-t$ and $\lambda\in[0,1]$, convexity also gives $\lambda K_s+(1-\lambda)K_{s'} \subseteq K_{\lambda s+(1-\lambda)s'}$.
Since $\gamma_D(K_s)=F(s)\in(0,1)$ on $[-t,\infty)$, Lemma~\ref{lem:ehrhard}, together with this inclusion and the monotonicity of $\Phi^{-1}$, yields
\begin{equation}
 \widetilde F(\lambda s+(1-\lambda)s') \ge\Phi^{-1}\bigl(\gamma_D(\lambda K_s+(1-\lambda)K_{s'})\bigr) \ge\lambda\widetilde F(s)+(1-\lambda)\widetilde F(s').
\end{equation}
Thus, we obtain that $\widetilde F$ is concave on $[-t,\infty)$.

Next, we define the anchor points $s_1\coloneqq\log(1-2\sigma)$ and $s_2\coloneqq\log(1+2\sigma)$ for the secant.
Since $\sigma\le1/4$, we have $s_1\ge-\log2>-1$, so that $-t<s_1<s_2$.
Moreover, since $2\sigma\le1/2$, the elementary bounds
$\log(1+2\sigma)\le2\sigma$ and
$-\log(1-2\sigma)\le4\sigma$ yield
\begin{equation}
 s_2-s_1=\log\frac{1+2\sigma}{1-2\sigma}\le6\sigma.
\end{equation}

We now bound $\widetilde F(s_1)$ and $\widetilde F(s_2)$ to control the secant slope.
Chebyshev's inequality gives
$\Pp\{\abs{e^{f(\mathbf G)}-1}\ge2\sigma\}\le1/4$, and hence
\begin{equation}
 F(s_1)=\Pp\{e^{f(\mathbf G)}\le1-2\sigma\}\le\frac14,\qquad
 1-F(s_2)=\Pp\{e^{f(\mathbf G)}>1+2\sigma\}\le\frac14.
\end{equation}
Applying $\Phi^{-1}$ to these bounds yields $\widetilde F(s_1)\le-q$ and $\widetilde F(s_2)\ge q$.
Together with $s_2-s_1\le6\sigma$, this shows that the secant through $(s_1,\widetilde F(s_1))$ and $(s_2,\widetilde F(s_2))$ has slope at least $q/(3\sigma)$.
To the left of $s_1$, concavity places the graph of $\widetilde F$ below the secant, see Figure~\ref{fig:secant}.
Consequently, using $\widetilde F(s_1)\le-q$ and $s_1>-1$, we obtain
\begin{equation}
 \widetilde F(-t)
 \le\widetilde F(s_1)
 -(s_1+t)\frac{\widetilde F(s_2)-\widetilde F(s_1)}{s_2-s_1}
 \le-\frac{q(t-1)}{3\sigma}.
\end{equation}
Applying $\Phi$ on both sides proves the first inequality in~\eqref{eq:convextail} by monotonicity.
The second one follows from $\Phi(-x)\le e^{-x^2/2}$ for $x\ge0$.
\end{proof}

\begin{remark}\label{rem:convextail}
The dependence on $\sigma^{-2}$ of the exponent
in~\eqref{eq:convextail} is sharp in order.
For $f_s(\mathbf g)=s g_1-s^2/2$, we have $\E[e^{f_s(\mathbf G)}]=1$ and $\sigma_s^2\coloneqq\Var(e^{f_s(\mathbf G)})=e^{s^2}-1\sim s^2$.
Since $f_s(\mathbf G)$ is Gaussian with variance $s^2$, for every fixed $t>1$, it holds
\begin{equation}
 -\log\Pp\{f_s(\mathbf G)\le-t\}
 \sim\frac{t^2}{2s^2}
 \sim\frac{t^2}{2\sigma_s^2}
 \qquad\text{as }s\downarrow0.
\end{equation}

Similar tail bounds have recently appeared in the literature.
Under the additional assumption $\E[e^{p f(\mathbf G)}]<\infty$ for every $p>0$, Chen~\cite{Chen26} gives a sharper lower-tail bound.
Paouris and Valettas~\cite[Thm.~2.1]{PV18} give a related bound in terms of the expected positive deviation above a median.
We provide this self-contained variant for the convenience of the reader.
\end{remark}

\subsection{Proof of Theorem~\ref{lem:smoothconfinement}}

\begin{proof}
Let $\kappa_{d,\psi}$ and $C_{d,\psi}$ be the constants from
Lemma~\ref{lem:expvariance}, and set
$a\coloneqq\frac14\sqrt{\kappa_{d,\psi}c_{d,\psi}}$ and
$\tau\coloneqq\kappa_{d,\psi}(n+1)$.
We need to show that the event $\{\norm{V_n}_\infty\le a\sqrt{nS_n}\}$ is exponentially rare.

First, we show that this event forces $\log U_\tau$ to be substantially below zero.
On this event, \eqref{eq:goodvariance} gives
$\abs{Y_{\mathbf x}}\le a\sqrt{n/c_{d,\psi}}$
for $\mathbf x\in A_n$.
Together with $\cosh u\le e^{\abs u}$ and the choices of $a$ and $\tau$, we obtain
\begin{equation}\label{eq:supnorm}
 \log U_\tau = -\frac{\tau}{2} + \log\biggl( \int_{A_n}\cosh(\sqrt\tau\,Y_{\mathbf x})\,d\nu_n(\mathbf x)\biggr)
 \le-\frac{\tau}{2}+a\sqrt{\frac{\tau n}{c_{d,\psi}}}
 \le-\frac{\tau}{4}.
\end{equation}
Thus, it holds $\Pp\{\norm{V_n}_\infty\le a\sqrt{nS_n}\} \leq \Pp\{ \log U_\tau
 \le-\frac{\tau}{4}\}$.
To bound $\Pp\{ \log U_\tau
 \le-\frac{\tau}{4}\}$, we apply Proposition~\ref{lem:convextail}.
By Lemma~\ref{lem:expvariance}, $\log U_\tau$ is a convex
function of $\mathbf G$, $\E[U_\tau]=1$, and
\begin{equation}
 \Var(U_\tau)\le\sigma^2
 \coloneqq C_{d,\psi}\frac{\tau^2}{\vert\D_n\vert}.
\end{equation}
Since $\tau^2/\vert\D_n\vert\to0$ and $\tau\to\infty$, we have $\sigma\le1/4$ and $\tau/4-1\ge\tau/8$ for sufficiently large $n$.
Thus, applying Proposition~\ref{lem:convextail} with $t\coloneqq\tau/4 \geq 1$ yields
\begin{equation}
 \Pp\Bigl\{\log U_\tau\le-\frac{\tau}{4}\Bigr\}
 \le\exp\Bigl(-\frac{q^2(t-1)^2}{18\sigma^2}\Bigr)
 \le\exp\Bigl(-\frac{q^2\vert\D_n\vert}{1152C_{d,\psi}}\Bigr).
\end{equation}
This proves~\eqref{eq:smoothconfinement} with
$b\coloneqq q^2/(1152C_{d,\psi})$.
\end{proof}

\begin{remark}\label{rem:generalatoms}
{\em (i)} Theorem~\ref{lem:smoothconfinement} extends to any real measurable atoms $a_R$ supported on $R\in\D_n$ with $\norm{a_R}_\infty\le1$ for which there exists $c_0>0$ independent of $n$ such that
\begin{equation}
 \sum_{R\in\D_n}\norm{a_R}_{L_2}^2\ge c_0S_n.
\end{equation}
The constants $a$ and $b$ depend only on $d$ and $c_0$.
Indeed, these assumptions imply $v_n\le S_n$ and $\bar v_n\ge c_0S_n$, so that the variance bound and the convexity argument apply unchanged.
Examples include $L_\infty$-normalized tensor Haar functions and nonnegative atoms such as hat functions.
Smoothness enters only in the subsequent function-space norm estimates.

{\em (ii)}  Deterministic small ball inequalities for Haar functions or atoms with vanishing moments must control cancellation between shapes for every coefficient vector.
Theorem~\ref{lem:smoothconfinement} instead shows that Gaussian coefficients produce a sum with small uniform norm only with exponentially small probability.
This probabilistic conclusion requires no vanishing-moment conditions.
\end{remark}

\begin{remark}[The case $d=1$]\label{rem:onedim}
For $d=1$, $V_n$ consists of $2^n$ independent bumps with disjoint supports.
Hence, it holds
\begin{equation}
 \Pp\{\norm{V_n}_\infty\le a\sqrt n\}
 \ge(1-e^{-a^2n/2})^{2^n}=e^{-o(\vert \D_n \vert)},
\end{equation}
so Theorem~\ref{lem:smoothconfinement} fails in this case.
The proof uses $d\ge2$ in the bound
$S_h/S_n\le(h+1)/(n+1)$ for $0\le h\le n$
from Lemma~\ref{lem:expvariance}.
\end{remark}
\section{Entropy and Kolmogorov numbers}\label{sec:entropy}\label{sec:lower}

Let $X,Y$ be Banach spaces and write $B_X\coloneqq \{x\in X:\norm x_X\le1\}$.
For a bounded set $K\subset Y$ and $t>0$, let $N(K,t;Y)$ be the smallest number of closed balls of radius $t$ with centers in $Y$ that cover $K$, possibly infinite.
Moreover, the Kolmogorov widths of $K$ are defined as
\begin{equation}
 d_k(K;Y)\coloneqq
 \inf_{\substack{V\subset Y\text{ linear}\\\dim V<k}}
 \sup_{f\in K}\inf_{v\in V}\norm{f-v}_Y.
\end{equation}
For a bounded linear operator $T\colon X\to Y$ and $k\in\N$, the dyadic entropy, Kolmogorov and approximation numbers are defined, respectively, by
\begin{equation}\label{eq:operatornumbers}
 \begin{aligned}
 e_k(T)
 &\coloneqq \inf\bigl\{t>0:N(T(B_X),t;Y)\le2^{k-1}\bigr\},\\
 d_k(T)
 &\coloneqq d_k(T(B_X);Y),\\
 a_k(T)
 &\coloneqq \inf_{\substack{A\colon X\to Y\text{ bounded linear}\\\operatorname{rank}A<k}}
       \norm{T-A}_{X\to Y}.
 \end{aligned}
\end{equation}
Carl's inequality~\cite{Carl81} is commonly used to derive lower bounds for approximation, Kolmogorov and Gelfand numbers from lower bounds for entropy numbers.
In the form of~\cite[eq.~(1.1)]{HKV16}, this inequality states that for every $\rho>0$, it holds
\begin{equation}\label{eq:carl}
 \max_{1\le k\le n}k^\rho e_k(T)
 \le C_\rho\max_{1\le k\le n}k^\rho d_k(T),
 \qquad n\in\N,
\end{equation}
where $C_\rho$ depends only on $\rho$.
Using this relation, we can establish the following lemma.
For the convenience of the reader, we give a self-contained proof.
The result is known in a more general situation,  see Cobos and K\"uhn \cite[Thm.\ 3.3]{CoKu09}.

\begin{lemma}\label{lem:carl}
Let $T\colon X\to Y$ be bounded, $\alpha>0$, $\beta\in\R$ and
$\eta_k\coloneqq k^{-\alpha}(\log(k+1))^\beta$.
If $d_k(T)\lesssim\eta_k$ and $e_k(T)\gtrsim\eta_k$ asymptotically, then it holds
\begin{equation}\label{eq:carlconsequence}
 d_k(T)\asymp e_k(T)\asymp\eta_k
 \qquad\text{as }k\to\infty.
\end{equation}
\end{lemma}

\begin{proof}
We use~\eqref{eq:carl} with $\rho=2\alpha$, see also~\cite[Lem.~4.1]{TU21}.
By assumption, there exist constants $C_1$, $c_1$ and $k_0 \in \N$ such that $d_k(T)\le C_1\eta_k$ and $e_k(T)\ge c_1\eta_k$ for all $k\ge k_0$.
After possibly enlarging $k_0$, the sequence $k^{2\alpha}\eta_k=k^\alpha(\log(k+1))^\beta$ is increasing for $k\ge k_0$.
On the other hand, it holds $k^{2\alpha}d_k(T)\le k_0^{2\alpha}\norm T$ for $k<k_0$.
Hence, for all $n \geq k_0$, we get
\begin{equation}\label{eq:carlweight}
 \max_{1\le k\le n}k^{2\alpha}d_k(T)\le C\,n^{\alpha}\bigl(\log(n+1)\bigr)^\beta.
\end{equation}

\noindent \emph{Upper bound for $e_k$.}
Using the bounds~\eqref{eq:carl} and~\eqref{eq:carlweight}, it holds that $n^{2\alpha}e_n(T)\le C\,n^\alpha(\log(n+1))^\beta$, namely that $e_n(T)\lesssim\eta_n$ for $n \geq k_0$.

\noindent\emph{Lower bound for $d_k$.}
The entropy lower bound $e_k(T)\ge c_1\eta_k$ and~\eqref{eq:carl} give
\begin{equation}\label{eq:carlchain}
 c_1\,n^{\alpha}(\log(n+1))^\beta
 \le n^{2\alpha}e_n(T)
 \le C_{2\alpha}\max_{1\le k\le n}k^{2\alpha}d_k(T).
\end{equation}
To lower bound a single Kolmogorov width, fix an integer $\lambda\ge2$ with $\lambda^\alpha\ge 2^{1+\abs\beta}C_{2\alpha}C/c_1$.
For every integer $m\ge\max\{\lambda,k_0\}$, put $n\coloneqq\lambda m$ and split the maximum at $m$.
Since $m\ge\lambda$, we have $n+1\le(m+1)^2$ and hence $\frac12\log(n+1)\le\log(m+1)\le\log(n+1)$.
The estimate~\eqref{eq:carlweight} and the choice of $\lambda$ therefore give
\begin{equation}
 \max_{k\le m}k^{2\alpha}d_k(T)
 \le C\,m^\alpha\bigl(\log(m+1)\bigr)^\beta
 \le \frac{c_1}{2C_{2\alpha}}
 n^\alpha\bigl(\log(n+1)\bigr)^\beta.
\end{equation}
Comparing this with~\eqref{eq:carlchain} shows that the maximum must be attained at some $m<k\le n$.
For these indices, monotonicity gives
$k^{2\alpha}d_k(T)\le n^{2\alpha}d_m(T)$, so that
\begin{equation}
 d_m(T)\ge\frac{c_1}{C_{2\alpha}}
 n^{-\alpha}\bigl(\log(n+1)\bigr)^\beta
 \ge\frac{c_1}{2^{\abs\beta}C_{2\alpha}\lambda^\alpha}\eta_m,
\end{equation}
since $n=\lambda m$ and $\log(m+1)\le\log(n+1)\le2\log(m+1)$.
\end{proof}

\subsection{From small balls to entropy}

For the covering argument, we need a probability bound for balls in
$L_\infty([0,1]^d)$ that is uniform in their centers.
Let $V_n'$ be an independent copy of $V_n$.
The coefficient differences are independent centered Gaussians of variance
$2$, so $V_n-V_n'$ has the same law as $\sqrt2\,V_n$.
For every $f\in L_\infty([0,1]^d)$ and $t>0$, independence and the triangle
inequality give
\begin{equation}\label{eq:shiftedball}
 \begin{aligned}
 \Pp\{\norm{V_n-f}_\infty\le t\}^2
 &=\Pp\bigl\{\norm{V_n-f}_\infty\le t,\ \norm{V_n'-f}_\infty\le t\bigr\}\\
 &\le\Pp\{\norm{V_n-V_n'}_\infty\le2t\}
 =\Pp\{\norm{V_n}_\infty\le\sqrt2\,t\}.
 \end{aligned}
\end{equation}
Set $\rho_n\coloneqq a\sqrt{nS_n}/\sqrt2$.
Combining~\eqref{eq:shiftedball} with Theorem~\ref{lem:smoothconfinement} gives $\Pp\{\norm{V_n-f}_\infty\le\rho_n\}\le e^{-b\vert\D_n\vert/2}$ for every center $f\in L_\infty([0,1]^d)$, namely that $V_n$ is uniformly unlikely to lie in any fixed ball.

Let $T\colon X\to Y$ be a bounded linear operator, where $X$ is a Banach space and $Y=C([0,1]^d)$ or $Y=L_\infty([0,1]^d)$.
Assume that its image ball $K\coloneqq T(B_X)$ is compact, and let $L_n>0$.
In the applications, $L_n$ will be chosen so that $V_n$ belongs to $L_n\cdot K$ with sufficiently large probability.
If $K$ is covered by $N$ balls of radius $\rho_n/L_n$ centered at $f_1,\ldots,f_N$, then scaling by $L_n$ and applying a union bound gives
\begin{equation}\label{eq:covering}
 \Pp\{V_n\in L_n\cdot K\}
 \le\sum_{j=1}^N\Pp\{\norm{V_n-L_nf_j}_\infty\le\rho_n\}
 \le N\,e^{-b\vert\D_n\vert/2}.
\end{equation}
For sufficiently large $n$, set $k_n\coloneqq\lfloor b\vert\D_n\vert/(4\log2)\rfloor$.
If $\Pp\{V_n\in L_n\cdot K\}\ge e^{-b\vert\D_n\vert/4}$, then~\eqref{eq:covering} forces $N\ge e^{b\vert\D_n\vert/4}\ge2^{k_n}>2^{k_n-1}$.
By the definition of entropy numbers, we therefore obtain
\begin{equation}\label{eq:entropytransfer}
 \Pp\{V_n\in L_n\cdot K\}\ge e^{-b\vert\D_n\vert/4}
 \quad\Longrightarrow\quad
 e_{k_n}(T)\ge\frac{\rho_n}{L_n}
 =\frac{a\sqrt{nS_n}}{\sqrt2\,L_n}.
\end{equation}
Here, $k_n\asymp2^nn^{d-1}$ and $k_{n+1}/k_n$ is bounded.
Monotonicity therefore extends the bounds of the form $k^{-\alpha}(\log k)^\beta$ proved below along $(k_n)$ to all sufficiently large $k$.
Below, we choose $L_n$ such that $\Pp\{V_n\in L_n\cdot K\}\ge e^{-b\vert\D_n\vert/4}$, which yields entropy lower bounds via \eqref{eq:entropytransfer}.

\begin{proof}[Proof of Theorem~\ref{thm:CM}]
The image ball $K_{d,p}\coloneqq T_{d,p}(B_{L_p})$ is uniformly bounded and equicontinuous by Hölder's inequality, and closed by reflexivity of $L_p$, hence compact by Arzelà--Ascoli.
Put $h_n\coloneqq\partial_1\cdots\partial_dV_n$.
Since the atoms in $V_n$ are smooth and vanish on the boundary, repeated integration gives $V_n=T_{d,p}h_n$.
At each point, $h_n$ is a centered Gaussian of variance at most $2^{2n}S_n\norm{\psi'}_\infty^{2d}$, since at most one differentiated atom of each shape contributes.
The Gaussian moment formula and Tonelli's theorem give \smash{$\E[\norm{h_n}_p^p] \le C_{d,p,\psi}\,2^{np}S_n^{p/2}$}.
Set $L_n\coloneqq(2C_{d,p,\psi})^{1/p}2^n\sqrt{S_n}$.
Since $V_n=T_{d,p}h_n$, Markov's inequality yields
\begin{equation}
 \Pp\{V_n\in L_n\cdot K_{d,p}\}
 \ge\Pp\{\norm{h_n}_p\le L_n\}
 \ge1-\frac{\E[\norm{h_n}_p^p]}{L_n^p}
 \ge\frac12.
\end{equation}
Applying~\eqref{eq:entropytransfer} and using $k_n\asymp2^nn^{d-1}$ and $\log k_n\asymp n$, we obtain
\begin{equation}
 e_{k_n}(T_{d,p})
 \gtrsim2^{-n}\sqrt n
 \asymp k_n^{-1}(\log k_n)^{d-1/2}.
\end{equation}

For $p\ge2$, the inclusion $L_p([0,1]^d)\hookrightarrow L_2([0,1]^d)$ and~\cite{DKLL99} give $d_k(T_{d,p})\lesssim k^{-1}(\log k)^{d-1/2}$.
Lemma~\ref{lem:carl} yields both asserted orders in this range.
For $1<p<2$, bounded mixed Sobolev extension and periodization map $K_{d,p}$ into a fixed multiple of the periodic mixed Sobolev unit ball.
Decomposing it into coordinatewise mean-zero components, applying Belinsky's entropy upper bound~\cite[Thm.~1.1]{Belinsky98} in each corresponding dimension, and covering the constant component gives the matching entropy upper bound.
\end{proof}

\subsection{Periodic spaces with mixed smoothness}\label{sec:periodic}\label{sec:spaces}

Choose an even $\varphi_0\in C_c^\infty(\R)$ that equals one on $[-1,1]$ and zero outside $[-2,2]$, and set $\varphi_s(\xi)\coloneqq \varphi_0(2^{-s}\xi)-\varphi_0(2^{1-s}\xi)$ for $s\ge1$.
Then, it holds $\sum_{s\ge0}\varphi_s(\xi)=1$.
For $f\in L_1(\mathbb T^d)$, we define
\begin{equation}\label{eq:dyadicblocks}
 \Delta_{\mathbf s}f(\mathbf x)
 \coloneqq \sum_{\mathbf k\in\mathbb Z^d}
       \biggl(\prod_{i=1}^d\varphi_{s_i}(k_i)\biggr)
       \widehat f(\mathbf k)e^{2\pi i\langle\mathbf k,\mathbf x\rangle}.
\end{equation}
Thus, $\Delta_{\mathbf s}$ selects frequencies of size $2^{s_i}$ in coordinate
$i$, with $s_i=0$ including the low frequencies.
Using these dyadic blocks, we define the norms
\begin{equation}\label{eq:mixednorms}
 \norm f_{\bW^r_p} \coloneqq \Biggl\|\biggl(\sum_{\mathbf s\in\N_0^d}
     2^{2r\abs{\mathbf s}}\bigl|\Delta_{\mathbf s}f\bigr|^2
     \biggr)^{1/2}\Biggr\|_p, \qquad \norm f_{\bB^r_{p,\theta}}\coloneqq \Biggl(\sum_{\mathbf s\in\N_0^d}
     2^{r\theta\abs{\mathbf s}}
     \bigl\|\Delta_{\mathbf s}f\bigr\|_p^\theta
     \Biggr)^{1/\theta}.
\end{equation}
Here, $\bW^r_p$ is the periodic mixed Sobolev (Bessel-potential) space for $1<p<\infty$ and $\bB^r_{p,\theta}$ is the mixed Besov space. For $\theta=\infty$, the mixed Besov norm is the supremum of the weighted block norms, and then $\bH^r_p\coloneqq \bB^r_{p,\infty}$ denotes the mixed H\"older-Nikol'skii space. The unit balls in $\bW^r_p$ are sometimes also called ``classes with bounded mixed derivative'', whereas the unit balls in $\bH^r_p$ are called ``classes with bounded mixed difference'', see \cite[Ch.\ 3]{DTU18}, \cite[Ch.\ 11]{TB04}. 
Other choices of $\varphi_0$ give equivalent norms.
For $r>1/p$, the canonical embeddings of these spaces into
$L_\infty(\mathbb T^d)$ are compact.
We denote them by $\operatorname{I}$.

In order to control the mixed smoothness norm of a hyperbolic layer we use the well-known atomic synthesis principle, see Vyb\' iral \cite[Thm.\ 2.4]{Vyb2006} and \cite[Prop.\ 5.5.10]{DTU18}.
The function space norm on the left-hand side of \eqref{eq:besovblocknorm} is controlled by a sequence space norm of the coefficients.

\begin{lemma}[Synthesis on a hyperbolic layer]\label{lem:synthesis}
Let $1<p\le\infty$, $1\le\theta\le\infty$ and $r>1/p$.
There is $C>0$, depending only on $d,p,\theta,r$ and $\psi$, such that for all $n\in\N_0$ and all real coefficients $(\lambda_R)_{R\in\D_n}$, it holds
\begin{equation}\label{eq:besovblocknorm}
 \biggl\|\sum_{R\in\D_n}\lambda_Ra_R\biggr\|_{\bB^r_{p,\theta}}
 \le C\,2^{rn}S_n^{1/\theta}\max_{R\in\D_n}\abs{\lambda_R}.
\end{equation}
\end{lemma}

\begin{proof}
For each shape $\mathbf s$, let $\D(\mathbf s)$ be its dyadic rectangles.
Apply the atomic synthesis estimate~\cite[Thm.~2.4(i)]{Vyb2006} on $\R^d$
to the atoms~\eqref{eq:atoms}, then use bounded periodization from their fixed
support to $\mathbb T^d$.
For $r>1/p$ and $p,\theta\ge1$, no moment conditions are required.
When $p=\infty$ or $\theta=\infty$, the corresponding $\ell_p$ or $\ell_\theta$ norm below is interpreted as a supremum.
This gives
\begin{equation}
    \biggl\|\sum_{R\in\D_n}\lambda_Ra_R\biggr\|_{\bB^r_{p,\theta}} \le C\,2^{(r-1/p)n} \biggl(\sum_{\abs{\mathbf s}=n} \biggl(\sum_{R\in\D(\mathbf s)}\abs{\lambda_R}^p\biggr)^{\theta/p} \biggr)^{1/\theta} \le C\,2^{rn}S_n^{1/\theta}\max_{R\in\D_n}\abs{\lambda_R}.
\end{equation}
The second inequality uses that there are $S_n$ shapes, each containing $2^n$ rectangles.
\end{proof}

\subsection{Proof of Theorem \ref{thm:periodicintro}}
We now combine Lemma~\ref{lem:synthesis} with the covering argument~\eqref{eq:entropytransfer} to obtain the entropy lower bounds in Theorem \ref{thm:periodicintro}.

\begin{proof}
Let $K$ be the unit ball of \smash{$\bB^r_{p,\theta}$}.
Choose a fixed $M>0$ with $\Pp\{\abs g\le M\}\ge e^{-b/4}$ for a standard Gaussian $g$.
Due to independence, we get
\begin{equation}\label{eq:coefficientbox}
 \Pp\Bigl\{\max_{R\in\D_n}\abs{G_R}\le M\Bigr\}
 =\bigl[\Pp\{\abs g\le M\}\bigr]^{\vert\D_n\vert}
 \ge e^{-b\vert\D_n\vert/4}.
\end{equation}
On this event, Lemma~\ref{lem:synthesis} gives
\begin{equation}
 \norm{V_n}_{\bB^r_{p,\theta}}
 \le L_n\coloneqq CM\,2^{rn}S_n^{1/\theta},
\end{equation}
so that $\Pp\{V_n\in L_n\cdot K\}\ge e^{-b\vert\D_n\vert/4}$.
Applying~\eqref{eq:entropytransfer} and using $S_n\asymp n^{d-1}$ and $k_n\asymp2^nn^{d-1}$, we obtain
\begin{equation}
 e_{k_n}\bigl(\operatorname I\colon\bB^r_{p,\theta}\to L_\infty(\mathbb T^d)\bigr)
 \gtrsim2^{-rn}\sqrt n\,S_n^{1/2-1/\theta}
 \asymp k_n^{-r}(\log k_n)^{(r+1/2-1/\theta)(d-1)+1/2}.
\end{equation}

For the mixed Sobolev lower bound, we use the embeddings
$\bB^r_{p,2}\hookrightarrow\bW^r_p$ for $2\le p<\infty$ and
$\bB^r_{2,2}=\bW^r_2\hookrightarrow\bW^r_p$ for $1<p<2$,
see~\cite[Ch.~3]{DTU18}. The unit ball of the smaller space is contained in
a multiple of the unit ball of $\bW^r_p$, so the mixed Besov lower bound with
$\theta=2$ applies.

For the upper bounds, we additionally require $r>1/2$, $p\ge2$ and $\theta\ge2$.
The Sobolev and mixed H\"older-Nikol'skii Kolmogorov upper bounds are proved in~\cite[Thms.~11.3.4 and~11.3.5]{TB04}.
For the Besov upper bound, write $I_q\colon\bB^r_{2,q}\to L_\infty(\mathbb T^d)$ for the embedding.
The endpoints $q=2$ and $q=\infty$ follow from $\bB^r_{2,2}=\bW^r_2$ and $\bB^r_{2,\infty}=\bH^r_2$.
For $2<\theta<\infty$, put $\eta\coloneqq1-2/\theta$.
The identity $(\bB^r_{2,2},\bB^r_{2,\infty})_{\eta,\theta}=\bB^r_{2,\theta}$, with equivalent norms, together with interpolation of Kolmogorov numbers~\cite[Thm.~3.2]{TU22}, yields
\begin{equation}
 d_{2k-1}(I_\theta)
 \lesssim d_k(I_2)^{2/\theta}d_k(I_\infty)^{1-2/\theta}
 \lesssim k^{-r}(\log k)^{(r+1/2-1/\theta)(d-1)+1/2}.
\end{equation}
The inclusion $\bB^r_{p,\theta}\hookrightarrow\bB^r_{2,\theta}$ gives the same upper bound for $p\ge2$.
Together with the derived entropy lower bounds, Lemma~\ref{lem:carl} yields the asserted equivalences.
\end{proof}

\begin{remark}\label{rem:ZeroMean}
The same proof applies to the coordinatewise mean-zero classes considered by Romanyuk~\cite{Romanyuk19}: choosing $\psi$ with $\int_0^1\psi(t)\,dt=0$ makes every atom mean-zero in each coordinate.
For $p\ge2$, the inclusion $\bB^r_{p,\theta}\hookrightarrow\bB^r_{2,\theta}$ transfers his entropy upper bound at $p=2$~\cite[Thm.~1]{Romanyuk19}.

For $1<p<2$ and $r>1/p$, Belinsky's upper bound~\cite[Thm.~1.1]{Belinsky98}, applied after the mean-zero decomposition from the proof of Theorem~\ref{thm:CM}, shows that the order for $e_k$ in Theorem~\ref{thm:periodicintro}(i) extends to this range.
This fails for $d_k$, which decays polynomially slower for $p<2$ \cite[Ch.~4]{DTU18}.
\end{remark}

\begin{remark}
Kolmogorov numbers can be substantially smaller than the errors of bounded linear approximation.
For $r>1/2$, let $\operatorname{I}_r\colon\bW^r_2\to L_\infty(\mathbb T^d)$ be the periodic mixed Sobolev embedding.
Its approximation numbers satisfy $a_k(\operatorname{I}_r) \asymp k^{-r+1/2}(\log k)^{r(d-1)}$, see \cite{CoKuSi16}.
Together with Theorem~\ref{thm:periodicintro}, this gives
\begin{equation}\label{eq:linearwidthgap}
 \frac{d_k(\operatorname{I}_r)}{a_k(\operatorname{I}_r)}
 \asymp\sqrt{\frac{\log k}{k}}.
\end{equation}
For Kolmogorov numbers, the approximant in a fixed subspace may depend nonlinearly on the input.
Approximation numbers require one bounded linear map of rank less than $k$.
\end{remark}

\section{From entropy to small deviations of Brownian sheets}\label{sec:sheet}

For a centered Gaussian random element $X$ in $C([0,1]^d)$, let $\mathcal H_X$ be the reproducing kernel Hilbert space associated with the covariance kernel $K_X(\mathbf s,\mathbf t)\coloneqq\E[X(\mathbf s)X(\mathbf t)]$.
Let $\operatorname{I}_X\colon\mathcal H_X\to C([0,1]^d)$ denote its Cameron--Martin embedding with $\operatorname{I}_Xh=h$, see~\cite[Sec.~2.4]{Bogachev98}.
The entropy numbers of $\operatorname{I}_X$ measure the size of the unit ball of $\mathcal H_X$ in the uniform norm.
For $r>1/2$ and $\beta\in\R$, the Kuelbs-Li-Linde correspondence \cite{KL93},  \cite[eq.~(1.3), Thm.~5.1]{LiLinde99}, gives
\begin{equation}\label{eq:entropybridge}
 e_k(\operatorname{I}_X)\asymp k^{-r}(\log k)^\beta  \Longleftrightarrow\quad
 -\log\Pp\{\norm X_\infty\le\varepsilon\}
 \asymp\varepsilon^{-1/(r-1/2)}
           \bigl(\log(1/\varepsilon)\bigr)^{\beta/(r-1/2)}.
\end{equation}
The forward implication also holds with $\gtrsim$ in place of $\asymp$ on both sides \cite[Thm.~1.1(III)]{LiLinde99}.

\begin{proof}[Proof of Theorem~\ref{thm:main}]
For the Brownian sheet, we get the covariance identity
\begin{equation}\label{eq:sheetCMcovariance}
 K_B(\mathbf s,\mathbf t) = \E[B(\mathbf s)B(\mathbf t)]
 =\langle\one_{[0,\mathbf s]},\one_{[0,\mathbf t]}\rangle_{L_2([0,1]^d)}.
\end{equation}
We now identify $\mathcal H_B$ with $T_d(L_2([0,1]^d))$.
Since $T_d$ is injective, its range becomes a Hilbert space $\mathcal H$ with inner product
\begin{equation}
 \langle T_dg,T_dh\rangle_{\mathcal H}
 =\langle g,h\rangle_{L_2}.
\end{equation}
The representation $(T_dg)(\mathbf t)=\langle g,\one_{[0,\mathbf t]}\rangle_{L_2}$ and~\eqref{eq:sheetCMcovariance} show that its reproducing kernel is indeed $K_B$.
Thus $\mathcal H=\mathcal H_B$ isometrically, with $\norm{T_dg}_{\mathcal H_B}=\norm g_{L_2}$.

Both $\operatorname{I}_B$ and $T_d$ map their unit balls onto the anchored mixed Sobolev ball $K_d$ in~\eqref{eq:Kd}.
Hence, Theorem~\ref{thm:CM} gives $e_k(\operatorname{I}_B)=e_k(T_d)\asymp k^{-1}(\log k)^{d-1/2}$.
Applying~\eqref{eq:entropybridge} with $r=1$ and $\beta=d-1/2$ proves~\eqref{eq:main}.
\end{proof}

\section{Further consequences}\label{sec:applications}

\subsection{Fractional integration and Brownian sheets}\label{sec:fractional}

Fractional integration provides a nonperiodic counterpart of the periodic spaces with dominating mixed smoothness considered in Subsection~\ref{sec:spaces}.
For $r>1/2$, we define the Riemann--Liouville operator by
\begin{equation}\label{eq:RLdefinition} 
 (R_rg)(t)\coloneqq\frac1{\Gamma(r)}\int_0^t(t-u)^{r-1}g(u)\,du
\end{equation}
and set $\smash{T_d^{(r)}}\coloneqq R_r^{\otimes d}\colon L_2([0,1]^d)\to C([0,1]^d)$.
This operator integrates with order $r$ in each coordinate.
Its image unit ball \smash{$K_d^{(r)}\coloneqq T_d^{(r)}(B_{L_2})$} models mixed smoothness $r$ with functions vanishing whenever any coordinate is zero.
More explicitly, it holds $(\smash{T_d^{(r)}}g)(\mathbf t)=\langle g,k_{\mathbf t}\rangle_{L_2}$, where $k_{\mathbf t}\in L_2([0,1]^d)$ is given by
\begin{equation}
 k_{\mathbf t}(\mathbf u)\coloneqq
 \begin{cases}
  \displaystyle\frac1{\Gamma(r)^d}\prod_{i=1}^d(t_i-u_i)^{r-1},
  &0\le u_i<t_i\text{ for every }i,\\
  0,&\text{otherwise}.
 \end{cases}
\end{equation}
The map $\mathbf t\mapsto k_{\mathbf t}$ is continuous in $L_2$.
For $\norm g_{L_2}\le1$, Cauchy--Schwarz gives $\norm{\smash{T_d^{(r)}}g}_\infty\le\sup_{\mathbf t}\norm{k_{\mathbf t}}_{L_2}<\infty$ and
\begin{equation}
 \bigl|(T_d^{(r)}g)(\mathbf t)-(T_d^{(r)}g)(\mathbf s)\bigr|
 \le\norm{k_{\mathbf t}-k_{\mathbf s}}_{L_2}.
\end{equation}
Therefore, the image unit ball is uniformly bounded and equicontinuous.
Since it is also closed by reflexivity of $L_2$, Arzelà--Ascoli shows that $\smash{K_d^{(r)}}$ is compact in the uniform norm.
For $r=1$, we recover $T_d$ and its image unit ball $K_d$.

\begin{theorem}\label{thm:fractionalentropy}
Let $d\ge2$.
For $r>1/2$, it holds
\begin{equation}\label{eq:fractionalentropy}
 e_k\bigl(T_d^{(r)}\bigr)\asymp d_k\bigl(T_d^{(r)}\bigr)
 \asymp k^{-r}(\log k)^{r(d-1)+1/2}
 \qquad\text{as }k\to\infty.
\end{equation}
The constants depend only on $r$ and $d$.
\end{theorem}

\begin{proof}
To apply~\eqref{eq:entropytransfer}, we construct the preimages of the atoms under $\smash{T_d^{(r)}}$ and estimate their $L_2$ norms.
Let $\psi$ be the fixed bump used in~\eqref{eq:atoms}, and choose an integer $q>r$.
Write $\psi_{s,m}\coloneqq\psi(2^s\,\cdot-m)$ and set $g_{r,s,m}\coloneqq R_{q-r}\smash{\psi_{s,m}^{(q)}}$.
Since $\smash{\psi_{s,m}^{(k)}}(0)=0$ for $k=0,\ldots,q-1$, the semigroup law $R_\alpha R_\beta=R_{\alpha+\beta}$ gives \smash{$R_rg_{r,s,m}=R_q\psi_{s,m}^{(q)}=\psi_{s,m}$}.
For $R=\prod_{i=1}^d I_{s_i,m_i}\in\D_n$, define $g_R(\mathbf u)\coloneqq\prod_{i=1}^d g_{r,s_i,m_i}(u_i)$.
Then $\smash{T_d^{(r)}}g_R=a_R$, and $g_n\coloneqq\sum_{R\in\D_n}G_Rg_R$ satisfies \smash{$T_d^{(r)}g_n=V_n$}.

To estimate $\E[\norm{g_n}_2^2]$, extend $g_{r,s,m}$ to $\R$ by replacing the lower limit $0$ in~\eqref{eq:RLdefinition} by $-\infty$, which preserves the values on $(0,1)$ since $\psi_{s,m}$ vanishes on $(-\infty,0]$.
This corresponds to the left-sided Liouville derivative of order $r$ of $\psi_{s,m}$, with Fourier multiplier of modulus $|\xi|^r$.
Writing $\mathcal F$ for the unitary Fourier transform, Plancherel's identity and scaling give
\begin{equation}\label{eq:fractionalatom}
 \norm{g_{r,s,m}}_{L_2(0,1)}^2
 \le\int_{\R}|\xi|^{2r}\bigl|(\mathcal F\psi_{s,m})(\xi)\bigr|^2\,d\xi
 =2^{(2r-1)s}\int_{\R}|\eta|^{2r}\bigl|(\mathcal F\psi)(\eta)\bigr|^2\,d\eta
 =C_{r,\psi}2^{(2r-1)s}.
\end{equation}
Since $\abs{\mathbf s}=n$, equation~\eqref{eq:fractionalatom} gives $\norm{g_R}_2^2\le C2^{(2r-1)n}$.
The Gaussian coefficients satisfy $\E[G_RG_{R'}]=0$ for $R\ne R'$ and $\E[G_R^2]=1$.
Therefore, the cross terms vanish when expanding $\norm{g_n}_2^2$ and taking the expectation.
Together with $\vert\D_n\vert=2^nS_n$, this yields
\begin{equation}\label{eq:fractionalenergy}
 \E\bigl[\norm{g_n}_2^2\bigr]
 =\sum_{R\in\D_n}\norm{g_R}_2^2
 \le C2^{2rn}S_n.
\end{equation}

We now apply~\eqref{eq:entropytransfer} to $\smash{K_d^{(r)}}$.
Set $L_n\coloneqq\sqrt{2C}\,2^{rn}\sqrt{S_n}$, where $C$ is the constant in~\eqref{eq:fractionalenergy}.
If $\norm{g_n}_2\le L_n$, then $V_n=\smash{T_d^{(r)}}g_n\in L_n\smash{K_d^{(r)}}$.
Thus, Markov's inequality and~\eqref{eq:fractionalenergy} give
\begin{equation}
 \Pp\{V_n\in L_nK_d^{(r)}\}
 \ge\Pp\{\norm{g_n}_2\le L_n\}
 \ge1-\frac{\E[\norm{g_n}_2^2]}{L_n^2}
 \ge\frac12.
\end{equation}
We can therefore apply~\eqref{eq:entropytransfer} for all sufficiently large $n$.
Using $k_n\asymp2^nn^{d-1}$ and $\log k_n\asymp n$, we obtain
\begin{equation}\label{eq:EntropyLower}
 e_{k_n}\bigl(T_d^{(r)}\bigr)
 \gtrsim2^{-rn}\sqrt n
 \asymp k_n^{-r}n^{r(d-1)+1/2}
 \asymp k_n^{-r}(\log k_n)^{r(d-1)+1/2}.
\end{equation}
This relation is extended to any $k$ as discussed after \eqref{eq:entropytransfer}.
Moreover, \cite[Thm.~5.7]{KuhnLinde02}, valid for every $r>1/2$, gives $\ell_n(T_d^{(r)})\lesssim n^{-r+1/2}(\log n)^{r(d-1)+1/2}$, and the inequality $d_{2n-1}(u)\lesssim n^{-1/2}\ell_n(u)$ for operators $u$ on Hilbert spaces~\cite{PTJ86} yields
\begin{equation}
 d_k\bigl(T_d^{(r)}\bigr)\lesssim k^{-r}(\log k)^{r(d-1)+1/2}.
\end{equation}
Consequently, Lemma~\ref{lem:carl} completes the proof.
\end{proof}

For $\mathbf H=(H_1,\ldots,H_d)\in(0,1)^d$, let $B^{\mathbf H}$ be the fractional Brownian sheet with covariance
\begin{equation}\label{eq:fractionalcovariance}
 \E\bigl[B^{\mathbf H}(\mathbf s)B^{\mathbf H}(\mathbf t)\bigr]
 =\prod_{i=1}^d\frac{s_i^{2H_i}+t_i^{2H_i}-\abs{s_i-t_i}^{2H_i}}2.
\end{equation}
Put $H_*\coloneqq\min_iH_i$ and $\nu\coloneqq\#\{i:H_i=H_*\}$.
When $H_i=H$ for all $i$, the Cameron--Martin unit ball of $B^{\mathbf H}$ contains $c\smash{K_d^{(H+1/2)}}$ for some $c>0$ by~\cite[Prop.~7.3 and the identity preceding~(7.9)]{KuhnLinde02}.
We also consider the integrated sheets $R_m^{\otimes d}B$, $m\in\N_0$, with $R_0$ being the identity.

\begin{theorem}\label{thm:fractionalsheet}
Let $d\ge2$, $\nu\ge2$ and $m\in\N_0$.
As $\varepsilon\downarrow0$, it holds
\begin{align}
 -\log\Pp\{\norm{B^{\mathbf H}}_\infty\le\varepsilon\}
 &\asymp\varepsilon^{-1/H_*}
       \bigl(\log(1/\varepsilon)\bigr)^{\nu-1+\nu/(2H_*)},\label{eq:fractionalsmallball}\\
 -\log\Pp\{\norm{R_m^{\otimes d}B}_\infty\le\varepsilon\}
 &\asymp\varepsilon^{-2/(2m+1)}
       \bigl(\log(1/\varepsilon)\bigr)^{d-1+d/(2m+1)}.\label{eq:integratedsmallball}
\end{align}
The constants depend only on $d$ and the respective regularity parameters.
\end{theorem}

\begin{proof}
The semigroup law gives $R_m^{\otimes d}T_d=\smash{T_d^{(m+1)}}$, so the Cameron--Martin unit ball of $R_m^{\otimes d}B$ is $\smash{K_d^{(m+1)}}$.
Theorem~\ref{thm:fractionalentropy} and~\eqref{eq:entropybridge}, with $r=m+1$ and $\beta=r(d-1)+1/2$, therefore prove~\eqref{eq:integratedsmallball}.

For the fractional sheet, first suppose that $H_1=\cdots=H_d=H$, and take $r=H+1/2$ and $\beta=r(d-1)+1/2$.
The unit-ball inclusion gives $e_k(\operatorname I_{B^{\mathbf H}})\ge c\,e_k(\smash{T_d^{(r)}})$.
Theorem~\ref{thm:fractionalentropy} and the lower-bound implication of~\eqref{eq:entropybridge} therefore give the lower estimate in~\eqref{eq:fractionalsmallball}.
For unequal Hurst parameters, fix all coordinates with $H_i>H_*$ at $1$, leaving a sheet with $\nu$ parameters and common Hurst parameter $H_*$.
Its supremum norm is bounded by that of the full sheet, so the preceding argument in dimension $\nu$, with $H=H_*$, gives the lower estimate in~\eqref{eq:fractionalsmallball}.
The matching upper estimate follows from~\cite[Cor.~7.5]{KuhnLinde02}, which covers both equal and unequal Hurst parameters.
\end{proof}

\subsection{Tensor-product Mat\'ern fields}

Tensor-product Mat\'ern kernels are used in Gaussian-process emulation and sparse-grid interpolation~\cite{AT26}.
Our estimates determine the entropy and Kolmogorov numbers of their Cameron--Martin embeddings and the corresponding uniform small-deviation rates.

\begin{corollary}\label{cor:matern}
Let $d\ge2$ and $X$ be a centered continuous Gaussian field on $[0,1]^d$ with covariance $\E[X(\mathbf s)X(\mathbf t)]=\prod_{i=1}^d k_i(s_i-t_i)$.
Suppose that, for a common $h>0$ and $\kappa_i>0$, the spectral density of $k_i$ is proportional to $\rho_i(\xi)=(\kappa_i^2+\xi^2)^{-h-1/2}$ for $\xi\in\R$.
Then, its Cameron--Martin embedding $\operatorname{I}_X\colon\mathcal H_X\to C([0,1]^d)$ satisfies
\begin{equation}\label{eq:maternentropy}
 e_k(\operatorname{I}_X)
 \asymp d_k(\operatorname{I}_X)
 \asymp k^{-h-1/2}(\log k)^{(h+1/2)(d-1)+1/2}
 \qquad\text{as }k\to\infty.
\end{equation}
Moreover, it holds
\begin{equation}\label{eq:maternsmallball}
 -\log\Pp\{\norm X_\infty\le\varepsilon\}
 \asymp\varepsilon^{-1/h}
 \bigl(\log(1/\varepsilon)\bigr)^{d-1+d/(2h)}
 \qquad\text{as }\varepsilon\downarrow0.
\end{equation}
The constants depend only on $d$, $h$ and the kernels $k_i$.
\end{corollary}

\begin{proof}
Put $r\coloneqq h+1/2$.
By~\cite[Prop.~1(b)]{AT26}, $\mathcal H_X$ is the mixed Sobolev space $\smash{\widehat{\bigotimes}_{i=1}^d}H^r(0,1)$, equipped with an equivalent norm.
Hence, we can obtain the entropy lower bound by comparing $\mathcal H_X$ with the periodic Sobolev space in Theorem~\ref{thm:periodicintro}.
For every $f \in \bW_2^r(\mathbb T^d)$, the restriction $f|_{[0,1]^d}$ lies in $\mathcal H_X$ and satisfies $\norm{f|_{[0,1]^d}}_{\mathcal H_X}\le C\norm f_{\bW_2^r}$.
Thus, multiplying the centers and radii of a uniform-norm cover of the $\mathcal H_X$ unit ball by $C$ gives a cover of the periodic Sobolev unit ball restricted to the cube.
The entropy lower bound in Theorem~\ref{thm:periodicintro} consequently yields
\begin{equation}
 e_k(\operatorname I_X)\gtrsim k^{-r}(\log k)^{r(d-1)+1/2}.
\end{equation}

For the upper bound, Sobolev extension in each coordinate, followed by cutoff and periodization, gives for each $f\in\mathcal H_X$ a periodic extension on a larger torus with Sobolev norm at most $C\norm f_{\mathcal H_X}$.
Rescaling the larger torus to $\mathbb T^d$ changes the Sobolev norms only by fixed factors.
Thus, Theorem~\ref{thm:periodicintro} provides subspaces $E_k\subset L_\infty(\mathbb T^d)$ of dimension less than $k$ that approximate the extended unit ball with error at most $Ck^{-r}(\log k)^{r(d-1)+1/2}$.
Now, we construct continuous subspaces with the same property. 
Convolution with a kernel from a smooth approximate identity maps $E_k$ linearly into $C(\mathbb T^d)$ and does not increase the dimension.
Moreover, these operators preserve the approximation bounds by contractivity on $L_\infty(\mathbb T^d)$.
Since they converge uniformly to the identity on the rescaled extensions of the $\mathcal H_X$ unit ball by relative compactness in $C(\mathbb T^d)$, the smoothing error can be made arbitrarily small.
Finally, undoing the rescaling and restricting the approximants to $[0,1)^d$ increases neither the dimension nor the uniform error, giving
\begin{equation}
 d_k(\operatorname I_X)\lesssim k^{-r}(\log k)^{r(d-1)+1/2}.
\end{equation}
Applying Lemma~\ref{lem:carl} proves~\eqref{eq:maternentropy}.
Since $r-1/2=h$ and $[r(d-1)+1/2]/h=d-1+d/(2h)$, the correspondence~\eqref{eq:entropybridge} gives~\eqref{eq:maternsmallball}.
\end{proof}

\subsection{Quantization of sheet laws}\label{sec:quantization}

Functional quantization approximates each realization of a random function by the closest function from a fixed finite collection of deterministic functions.
More precisely, for a continuous Gaussian field $X$ with law $\mu$ on $C([0,1]^d)$, $m\in\N$ and $p=1,2$, we define
\begin{equation}\label{eq:quantizationdefinition}
 q_{m,p}(\mu)\coloneqq
 \inf_{f_1,\ldots,f_m}
 \Bigl(\E\Bigl[\min_{1\le j\le m}\norm{X-f_j}_\infty^p\Bigr]\Bigr)^{1/p}
 =\inf_{\#\operatorname{supp}\eta\le m}W_p(\mu,\eta).
\end{equation}
The first infimum runs over $f_1,\ldots,f_m\in C([0,1]^d)$, the second over probability measures $\eta$, and $W_p$ is the $p$-Wasserstein distance induced by the uniform norm.
Small deviations determine quantization rates for broad classes of Gaussian measures~\cite{DFMS03,GLP03}.
Using~\cite[Thm.~1.4]{BLG14}, we compare optimal quantization with approximation by empirical measures.

\begin{corollary}\label{cor:sheetquantization}
Let $d\ge2$ and let $X$ be the Brownian sheet $B$, the fractional sheet $B^{\mathbf H}$ with $\mathbf H\in(0,1)^d$ and $\nu\ge2$, or the integrated sheet $R_q^{\otimes d}B$ with $q\in\N_0$.
Let $\mu$ be its law, let $X^{(1)},\ldots,X^{(m)}$ be independent copies, and put $\mu_m\coloneqq m^{-1}\sum_{i=1}^m\delta_{X^{(i)}}$.
Then, as $m\to\infty$, it holds
\begin{equation}\label{eq:empiricaltransport}
 q_{m,1}(\mu)\asymp q_{m,2}(\mu)\asymp\E[W_2(\mu,\mu_m)]\asymp
 \begin{cases}
  \displaystyle\frac{(\log\log m)^{d-1/2}}{(\log m)^{1/2}},
  &X=B,\\[6pt]
  \displaystyle\frac{(\log\log m)^{H_*(\nu-1)+\nu/2}}{(\log m)^{H_*}},
  &X=B^{\mathbf H},\\[6pt]
  \displaystyle\frac{(\log\log m)^{(q+1)(d-1)+1/2}}{(\log m)^{q+1/2}},
  &X=R_q^{\otimes d}B.
 \end{cases}
\end{equation}
The constants depend only on $d$ and the respective regularity parameters.
\end{corollary}

\begin{proof}
Since $\mu_m$ has at most $m$ support points, equation~\eqref{eq:quantizationdefinition} gives $q_{m,1}(\mu)\le q_{m,2}(\mu)\le W_2(\mu,\mu_m)$ for every realization.
Therefore, it suffices to bound $\E[W_2(\mu,\mu_m)]$ from above and $q_{m,1}(\mu)$ from below.
The small-deviation bounds~\eqref{eq:main}, \eqref{eq:fractionalsmallball} and~\eqref{eq:integratedsmallball} satisfy the doubling and inverse-growth hypotheses of~\cite[Thm.~1.4]{BLG14}.
For a small-deviation order $t^{-\alpha}(\log(1/t))^\beta$ with $\alpha>0$, that theorem gives $\E[W_2(\mu,\mu_m)]\lesssim(\log m)^{-1/\alpha}(\log\log m)^{\beta/\alpha}$, yielding the upper rates in~\eqref{eq:empiricaltransport}.
For $X=B^{\mathbf H}$, we use $\alpha=1/H_*$ and $\beta=\nu-1+\nu/(2H_*)$.

For the lower bound, the independent-copy argument in~\eqref{eq:shiftedball} gives
\begin{equation}
 \Pp\{\norm{X-f}_\infty\le t\}^2
 \le\Pp\{\norm X_\infty\le\sqrt2\,t\}
 \qquad\text{for }f\in C([0,1]^d),\ t>0.
\end{equation}
Now, set $t_m\coloneqq c(\log m)^{-1/\alpha}(\log\log m)^{\beta/\alpha}$.
The small-deviation estimates give $-\log\Pp\{\norm X_\infty\le\sqrt2\,t_m\}\gtrsim c^{-\alpha}\log m$, so a sufficiently small fixed $c>0$ ensures $\Pp\{\norm X_\infty\le\sqrt2\,t_m\}\le(2m)^{-2}$ for all sufficiently large $m$.
Every radius-$t_m$ ball then has mass at most $1/(2m)$, so the union of any $m$ such balls has mass at most $1/2$.
Consequently, $\E[\min_{1\le j\le m}\norm{X-f_j}_\infty]\ge t_m/2$ for every choice of representatives.
Taking the infimum proves the lower bound for $q_{m,1}(\mu)$.
\end{proof}
The Brownian-sheet rate extends the two-dimensional result of~\cite[p.~544]{BLG14}.
For fractional sheets with $\nu\ge3$, the quantization rates attain the order conjectured in~\cite[Sec.~3.1, after eq.~(3.6)]{LP07}.
The same rates give the optimal worst-case error for deterministic integration of $1$-Lipschitz functionals against these sheet laws, when only exact functional evaluations count towards the cost~\cite[Thm.~1]{CDMR09}.
\vspace{0.2cm}\\
\noindent\textbf{Acknowledgment.} 
DB is supported by Simons Travel Support for Mathematicians. TU is supported by the German Research Foundation (DFG) with grant Ul403/4-1. The authors would like to thank Thomas K\"uhn for further references and valuable remarks. Several of the presented relations were developed with the assistance of LLMs by Anthropic and OpenAI. Additional verification, presentation and writing were carefully carried out by the authors, who take full responsibility for the results.

\bibliographystyle{abbrvurl}
\bibliography{references}
\end{document}